\documentclass[3p,times]{elsarticle}

\usepackage[utf8]{inputenc}
\usepackage[american]{babel}
\usepackage{empheq}

\usepackage{pgfplots}
\pgfplotsset{compat=newest}
\pgfplotsset{/pgf/number format/1000 sep={\,}}
\definecolor{red}{RGB}{150,0,0}
\definecolor{green}{RGB}{0,150,0}
\definecolor{blue}{RGB}{0,0,150}
\definecolor{fillblue}{RGB}{110,150,121}
\newcommand{\opacityI}{0.1}
\newcommand{\opacityII}{0.5}
\newcommand{\opacityIII}{0.9}
\newcommand{\quantilelabels}[1]{%
\draw (5.7cm,3cm) node[]{$95\,\%$};
\draw[fill=#1, opacity=\opacityI] (4.8cm, 2.8cm) rectangle (5.2cm, 3.2cm);
\draw[fill=#1, opacity=\opacityII] (4.8cm, 2.8cm) rectangle (5.2cm, 3.2cm);
\draw[fill=#1, opacity=\opacityIII] (4.8cm, 2.8cm) rectangle (5.2cm, 3.2cm);
\draw[] (4.8cm, 2.8cm) rectangle (5.2cm, 3.2cm);
\draw (5.7cm,2.5cm) node[]{$99\,\%$};
\draw[fill=#1, opacity=\opacityI] (4.8cm, 2.3cm) rectangle (5.2cm, 2.7cm);
\draw[fill=#1, opacity=\opacityII] (4.8cm, 2.3cm) rectangle (5.2cm, 2.7cm);
\draw[] (4.8cm, 2.3cm) rectangle (5.2cm, 2.7cm);
\draw (5.7cm,2cm) node[]{$100\,\%$};
\draw[fill=#1, opacity=\opacityI] (4.8cm, 1.8cm) rectangle (5.2cm, 2.2cm);
\draw[] (4.8cm, 1.8cm) rectangle (5.2cm, 2.2cm);
}

\usepackage{amsmath,amssymb,amsthm,amsfonts}

\usepackage{colortbl}
\usepackage{mathtools}

\newcommand{\transp}{\ensuremath{\top}}
\usepackage{xspace}

\newcommand{\eg}{e.\,g.\ }
\newcommand{\ie}{i.\,e.\ }
\newcommand{\cf}{cf.\ }

\newcommand{\bs}[1]{\ensuremath{\boldsymbol{#1}}}
\newcommand{\obs}[1]{\ensuremath{\overline{\boldsymbol{#1}}}}
\newcommand{\tbs}[1]{\ensuremath{\widetilde{\boldsymbol{#1}}}}

\renewcommand{\bar}[1]{\ensuremath{\overline{#1}}}
\newcommand{\pfrac}[2]{\ensuremath{\frac{\partial #1}{\partial #2}}}
\newcommand{\lint}{\ensuremath{\int\limits}}
\newcommand{\Expect}[1]{\ensuremath{\mathbb{E}\left[#1\right]}}

\newcommand{\Variance}[1]{\ensuremath{\mathbb{V}\left[#1\right]}}
\newcommand{\Prob}[1]{\ensuremath{\mathbb{P}\left[#1\right]}}
\newcommand{\dummybox}{\fbox{\vbox to 4cm {\vfil \hbox to 5cm{\quad}\vfil}}}
\DeclarePairedDelimiter\floor{\lfloor}{\rfloor}
\begin{document}

\begin{frontmatter}




\title{Stochastic Nonlinear Model Predictive Control with\\Efficient Sample Approximation of Chance Constraints}

\author[TUI,IFAT,cor1]{Stefan Streif}
\author[IFAT]{Matthias Karl}
\author[BER]{Ali Mesbah}

\address[TUI]{Institute for Automation and Systems Engineering, Ilmenau University of Technology, 98684 Ilmenau, Germany}

\address[IFAT]{Institute for Automation Engineering, Otto-von-Guericke Universität Magdeburg, 39106 Magdeburg, Germany}

\address[BER]{Department of Chemical and Biomolecular Engineering, University of California, Berkeley, California 94720, USA}

\address[cor1]{Corresponding author (Stefan.Streif@TU-Ilmenau.de).}

\begin{abstract} 
This paper presents a stochastic model predictive control approach for nonlinear systems subject to time-invariant probabilistic uncertainties in model parameters and initial conditions. The stochastic optimal control problem entails a cost function in terms of expected values and higher moments of the states, and chance constraints that ensure probabilistic constraint satisfaction. The generalized polynomial chaos framework is used to propagate the time-invariant stochastic uncertainties through the nonlinear system dynamics, and to efficiently sample from the probability densities of the states to approximate the satisfaction probability of the chance constraints.
To increase computational efficiency by avoiding excessive sampling, a statistical analysis is proposed to systematically determine a-priori the least conservative constraint tightening required at a given sample size to guarantee a desired feasibility probability of the sample-approximated chance constraint optimization problem. 
In addition, a method is presented for sample-based approximation of the analytic gradients of the chance constraints, which increases the optimization efficiency significantly. 
The proposed stochastic nonlinear model predictive control approach is applicable to a broad class of nonlinear systems with the sufficient condition that each term is analytic with respect to the states, and separable with respect to the inputs, states and parameters. 
The closed-loop performance of the proposed approach is evaluated using the Williams-Otto reactor with seven states, and ten uncertain parameters and initial conditions.
The results demonstrate the efficiency of the approach for real-time stochastic model predictive control and its capability to systematically account for probabilistic uncertainties in contrast to a nonlinear model predictive control approaches.
\end{abstract}
%
%
%

\end{frontmatter}


\section{Introduction}
\label{sec:intro}

Model predictive control (MPC) is widely used in the process industry owing to its ability to deal with multivariable complex dynamics and to incorporate system constraints into the optimal control problem \cite{Qin_etAl_1997_Overview_industrial_MPC,Rawlings_2000_Tutorial_Overview_MPC}. However, parametric uncertainties and exogenous disturbances are ubiquitous in real-world systems, and the classical MPC framework is inherently limited to systematically account for uncertainties \cite{mor99}. This consideration has led to the development of numerous robust MPC formulations that deal with uncertainties. The robust MPC approaches can be broadly categorized as deterministic and stochastic approaches based on the representation of uncertainties and the handling of constraints.

In deterministic robust MPC approaches (for a review see, e.\,g., \cite{Bemporad_Morari_1999_RMPC_Survey}), uncertainties are often assumed to be bounded. The control law is determined such that the control objective is minimized with respect to worst-case uncertainty realizations, and/or such that the constraints are satisfied for all admissible values of uncertainties. Hence, robust MPC approaches discard statistical properties of uncertainties and are conservative \cite{Garatti_Campi_2013_CSM__Modulating_Robustness_ControlDesign,Vidyasagar_2001_Automatica__RA_Rob_ContrSyn} if the worst-case uncertainty realizations have a small probability of occurrence.

In stochastic MPC (SMPC) approaches (e.\,g., see early work \cite{sch99,Li_etAl_2000_RMPC_CC,hes03}) uncertainties are described by probability distributions (instead of bounded sets), which can often be readily obtained from state or parameter estimations.
Such a stochastic approach to MPC not only alleviates the conservatism of worst-case control, but also enables tuning robustness against performance by allowing prespecified levels of risk during operation. The trade-off between control performance and robustness is achieved using chance (or probabilistic) constraints, which ensure the satisfaction of constraints with a desired probability level. 

Stochastic MPC has recently become an active research area \cite{Li_etAl_2008_CC_PC_UncSys,Li_etAl_2000_RMPC_CC,Cannon_etAl_2011_Stochastic_TubeMPC_CC,pri09,ber09,Cannon_etAl_2009_Automatica__MPC_StochMultUnc,Zhang_etAl_2013_CDC__LinSMPC_RA_RO,Blackmore_etAl_2010_IEEETRob_Particle_SMPC_CC,Oldewurtel_etAl_2008_CDC_Approx_CC_affineDisturbFB,Kouvaritakis_etAl_2013_IJSS_DisturbComp_SMPC,Korda_etAl_2014_IEEETAC__LinSMPC_AverageConstrViol,Hashimoto_2013_CDC__LinSMPC_Chebychev,Farina_etAl_2013_CDC__LinSMPC_Cantelli}. These stochastic optimal control approaches are often limited to linear systems or restricted to certain types of uncertainty distributions (e.\,g.\ Gaussian uncertainties). The reason is twofold: first, the complexity of predicting the propagation of general uncertainty distributions through nonlinear system dynamics; second, chance constrained optimization problems are in general computationally intractable due to computation of multi-dimensional probability integrals. 
These integrals can only be evaluated exactly if special probability distributions are assumed (see, e.\,g., \cite{Calafiore_Ghaoui_20016_JOTA__DistrRobust_CC_LP,Geletu_etAl_2013_IJSS__Review_CCOpt}). 
However, such assumptions may not hold for the probability distributions of the states, especially in the presence of nonlinear system dynamics.

The restriction to special distributions is alleviated in sampling or scenario-based methods, or randomized algorithms (see, e.\,g.\ \cite{
Li_etAl_2000_RMPC_CC,
Royset_Polak_2004_SAA_CC,
Vidyasagar_2001_Automatica__RA_Rob_ContrSyn,
Shapiro_2008_MP__StochProgr_UncOpt,
Garatti_Campi_2013_CSM__Modulating_Robustness_ControlDesign,
Calafiore_Campi_2006_IEEETAC__Scenario_Approach_RobustControl,
Campi_Garatti_2008_SIAMJO__ExactFeasibility_RA_ConvexPrograms,
Calafiore_2010_SIAMJO__Random_ConvexPrograms,
Schildbach_etAl_2012_ACC__RA_LinMPC,
Calafiore_Fabiano_2013_IEEETAC__RMPC_ScenarioOpt,
Zhang_etAl_2013_CDC__LinSMPC_RA_RO}).
In scenario-based approaches, for instance, a suitable number of randomly extracted scenarios of uncertainty and disturbances are used to formulate an optimization problem that replaces the chance constrained control problem \cite{Garatti_Campi_2013_CSM__Modulating_Robustness_ControlDesign,Calafiore_Campi_2006_IEEETAC__Scenario_Approach_RobustControl,
Campi_Garatti_2008_SIAMJO__ExactFeasibility_RA_ConvexPrograms,
Calafiore_2010_SIAMJO__Random_ConvexPrograms,
Schildbach_etAl_2012_ACC__RA_LinMPC,
Calafiore_Fabiano_2013_IEEETAC__RMPC_ScenarioOpt,
Zhang_etAl_2013_CDC__LinSMPC_RA_RO}.
However, these approaches usually assume convexity of the optimization problem, which often implies linearity of the system dynamics and constraints.
Furthermore, high-dimensional uncertainties often require a large number of samples for accurate prediction of the system dynamics, which may be prohibitively expensive for real-time control.
Bounds on the required sample size in the scenario approach are given e.\,g.\ in \cite{Alamo_etAl_2010_Sample_Complexity_Prob_Analysis_Design,Alamo_etAl_2014_RA_UncSys_Sample_Complexity}.
Thus, efficient sampling and evaluation of the randomized constraints is crucial for the efficiency of these approaches.
\footnote{There exists a large literature on alternatives to sampling-based methods
\cite{Geletu_etAl_2014_EOpt__AnalyticApprox_NonConvex_CC,
Hashimoto_2013_CDC__LinSMPC_Chebychev,
Feng_etAl_2011_CDC_KinshipFuncApprox_CC,
Mesbah_etAl_2014_ACC__SNMPC_CC,
Farina_etAl_2013_CDC__LinSMPC_Cantelli}.
See also \cite{Geletu_etAl_2013_IJSS__Review_CCOpt} for a recent review.}

In addition to the need for accurate and computationally efficient approximation of chance constraints, the evaluation of gradients of chance constraints is also critical for real-time control applications. As shown \eg in \cite{Imsland_etAl_2010_MB_Contr_Comparison_FiniteDiff_SensODE}, providing analytic forms of the gradients of the objective function and constraints can significantly improve the speed and accuracy of online optimization. However, the computation of gradients for nonlinear chance constraints is particularly challenging for general probability distributions, as analytic expressions of the gradients cannot be readily derived \cite{Geletu_etAl_2013_IJSS__Review_CCOpt}.
For sample-based approaches, it is known that a finite-difference approximation of the gradients can be very inaccurate even for a large number of samples \cite{Garnier_etAl_2009_Approximative_Gradients_CC,Royset_Polak_2004_SAA_CC}.

The contribution of this work is a stochastic nonlinear MPC (SNMPC) framework based on sample approximation of the chance constraints (Sec.~\ref{sec:CC}) and their gradients (Sec.~\ref{sec:gradCC}).
In addition, a statistical analysis is presented to determine a-priori (i.e., before the real-time optimizations) the required constraint tightening and the number of samples that guarantee a desired feasibility probability for a prespecified robustness (or risk) level (Sec.~\ref{sec:CC}).
The presented SNMPC approach is applicable to a broad class of nonlinear systems subject to time-invariant uncertainties in model parameters and initial conditions. The system dynamics and constraints are required to be analytic with respect to the states and separable with respect to the inputs, states and parameters.
The generalized polynomial chaos (PC) framework is used to obtain a computationally efficient surrogate for uncertainty propagation through the nonlinear system dynamics in order to generate a large number of Monte-Carlo (MC) samples.%
\footnote{The use of the PC framework for stochastic MPC and optimal control has also been investigated in \cite{Mesbah_etAl_2014_ACC__SNMPC_CC,Paulson_etAl_2014_Fast_SMPC,fag12,Huschto_Sager_2013_ECC__OC_PCE,Kwang-Ki_Braatz_2013_IJC_MPC_PCE}.}
In the PC framework, spectral expansions in terms of orthogonal polynomials are used to present stochastic quantities \cite{wie38,Ghanem_Spanos_1991_SFE,xiu02,Kim_etal_2013_CSM__PCE_review}, which allows sampling in a computationally efficient manner (Sec.~\ref{sec:PC}).
The performance of the proposed SNMPC approach is demonstrated for the Williams-Otto reactor using extensive simulation studies (Sec.~\ref{sec:example}).


\subsubsection*{Notation}

Bold symbols (e.\,g.\ $\bs{x}$, $\bs{\xi}$) denote vectors.
Subscripts at vectors are used for indexing vector elements.
$n_x$, $n_u$, \ldots\ denote dimension of the vectors indicated by the subscripts.
$N$ (e.\,g.\ $N_x$, $N_S$) denote the number of constraints or samples. 
Sets are denoted by calligraphic letters (e.\,g.\ $\mathcal{X}$,  $\mathcal{U}$).
Tilde (e.\,g.\ $\widetilde{\bs x}$, $\widetilde{\bs p}$) denotes coefficients in the polynomial chaos expansion of the corresponding variable, and $\widetilde{P}$ denotes the number of terms/coefficients in the polynomial chaos expansion of order $P$.
Probability distributions are denoted by $\mu$.
Expectation of a random variable is denoted by $\Expect{\cdot}$, variance by $\Variance{\cdot}$, higher order moments by $\mathbb{E}^m[\cdot]$, and probability by $\Prob{\cdot}$. 
Superscripts $^{[i]}$ (e.\,g.\ $\bs{\xi}^{[i]}$) indicate independent and identically distributed samples.
Superscript $^\ast$ denotes the optimal solution of an optimization problem.

\section{Stochastic Model Predictive Control Problem}
\label{sec:problem}

Consider the continuous-time, nonlinear system
\begin{align}
\dot{\bs{x}}(t) &= \bs{f}(\bs{x}(t),\bs{u}(t),\bs{p}), \; \bs{x}(0) = \bs{x}_0 \label{eq:sys:f},
\end{align}
where $t$, $\bs{x} \in \mathbb{R}^{n_x}$, $\bs{u} \in \mathbb{R}^{n_u}$ and $\bs{p}\in\mathbb{R}^{n_p}$ denote time, the states, the inputs and the time-invariant parameters, respectively. $\bs{x}_0$ denotes the initial states. The function $\bs{f} : \mathbb{R}^{n_x}\times\mathbb{R}^{n_u}\times\mathbb{R}^{n_p} \rightarrow \mathbb{R}^{n_x}$ represents the nonlinear system dynamics.
To be able to efficiently employ the framework presented in the subsequent sections, it is assumed that $\bs{f}$ can be transformed into a polynomial-in-the-states representation \cite{Ohtsuka_2005_IEEETAC_Immersion}. A sufficient condition for this assumption to hold is that $\bs{f}$ is analytic with respect to the states, and separable with respect to the inputs, states and parameters.

\subsection{Uncertainties} 

The system dynamics are subject to the following uncertainties. The time-invariant parameters $p_i$, $i=1,\ldots, n_p$ are assumed to be distributed according to the continuous probability distributions $\mu_{p_i}$ (denoted by $p_i \sim \mu_{p_i}$). Additionally, uncertain estimates of the states $x_i(t_k)$, $i=1,\ldots, n_x$ described by continuous probability distributions $\mu_{x_i(t_k)}$, can be used to recursively initialize \eqref{eq:sys:f}. Such uncertainty descriptions for the initial conditions are often available from state estimation, for example from Kalman filters. Note that this formulation also allows considering exact state estimates by choosing $\mu_{x_i(t_k)}$ to be a Dirac distribution. For technical reasons, it is assumed that the parameters and initial conditions are uncorrelated and have finite variances (e.\,g., $\Variance{p_i}< \infty$ and $\Variance{x_i(t_k)}< \infty$).

\subsection{Cost Function and Constraints} 

This work considers the stochastic optimal control of the system \eqref{eq:sys:f} on the finite-time horizon $[t_k,t_f]$
\footnote{\label{fn:horizon}Here $t_f = t_k + T$ (with $T$ being the prediction horizon) for receding horizon control, and fixed $t_f$ for shrinking horizon control.}%
, while constraints on the inputs and states should be satisfied in the presence of uncertainties.
The cost function of the stochastic nonlinear model predictive control approach is assumed to be deterministic and defined by
\begin{equation}
\label{eq:cost}
J(\bs{x}(\cdot), \bs{u}(\cdot)) 
\coloneqq \\
\int_{t_k}^{t_f} F\left(\bs{u}(t), \mathbb{E}[\bs{x}(t)], \mathbb{E}^2[\bs{x}(t)], \ldots \right) \text{d}t \\
+ E\left(\mathbb{E}[\bs{x}(t_f)], \mathbb{E}^2[\bs{x}(t_f)], \ldots \right).
\end{equation}
Here $F$ and $E$ denote the running and terminal cost functions, respectively, both of which can be functions of the moments of the states $\mathbb{E}^{m}[\cdot]$. 
Such a cost function enables shaping state distributions or, in a simpler case, minimizing the variance of state distributions (e.g., see \cite{Mesbah_etAl_2014_ACC__SNMPC_CC,Fisher_Bhattacharya_2009_LQR_PCE}).

In the following, we assume that \eqref{eq:sys:f} and \eqref{eq:cost} are time-discretized to integrate the nonlinear system dynamics and to impose constraints on the states and inputs at discrete time-points as described next.

State constraints\footnote{Output constraints can be considered similarly.} are imposed at $N_x$ different time points $t_{x,i} \in [t_k, t_f]$, $i = 1, \ldots, N_{x}$
\begin{align}
\label{eq:cons:x}
\bs{X} \in \mathcal{X} \coloneqq \{ 
g_i(\bs{X}) \le 0,\quad i = 1, \ldots, N_g 
 \} \subset \mathbb{R}^{n_x\,N_x},
\end{align}
where $\bs{X} \coloneqq \left[ \bs{x}(t_{x,1})^\transp, \ldots, \bs{x}(t_{x,N_x})^\transp  \right]^\transp$,
and $N_g$ is the number of constraints.
To efficiently employ the proposed control approach, functions $g_i(\bs X)$ are assumed to satisfy the same conditions as $\bs{f}$ (i.e., being analytic with respect to the states and separable with respect to the inputs, states and parameters).

In addition to state constraints, inputs are constrained by a compact set $\mathcal{U}$
\begin{align}
\label{eq:cons:u}
\bs u(t) \in \mathcal{U}\subset\mathbb{R}^{n_u},\quad \forall t \in [t_k, t_f].
\end{align}
For notational simplicity, output constraints are not explicitly considered here, as output constraints can often be represented in terms of \eqref{eq:cons:x}. 
Note that algebraic equations can be straightforwardly incorporated into the considered stochastic optimal control framework \cite{Paulson_etAl_2014_Fast_SMPC}.

Under the uncertainties in the system parameters and initial conditions, the solution trajectories of system \eqref{eq:sys:f} may violate the constraints \eqref{eq:cons:x}.
In this work, inputs $\bs{u}(t)$ are designed to satisfy \eqref{eq:cons:u} such that constraints \eqref{eq:cons:x} are fulfilled in a probabilistic manner in the presence of uncertainties.
This is formalized by chance constraints \cite{Li_etAl_2008_CC_PC_UncSys}
\begin{align}
\label{eq:cons:CC}
\Prob{\bs{g}(\bs{X}) \le \bs{0}} \ge \beta
\end{align}
where $\beta \in (0,1]$ is a user-specified probability chosen according to process requirements: $\beta=1$ corresponds to hard constraints that should hold at all times for all uncertainty realizations (i.e., the risk-free case); $\beta<1$ allows for constraint violation with probability $1-\beta$ in order to trade-off control performance with robustness.

When $N_g > 1$, \eqref{eq:cons:CC} entails \textit{joint chance constraints}, as all constraints $g_1(\bs X) \le 0, \ldots, g_{N_g}(\bs X) \le 0$ should be jointly satisfied with the probability level $\beta$. If $N_g = 1$ or constraints are defined independently using different $\beta_i$ for each constraint $g_i(\bs X)$, \eqref{eq:cons:CC} will be referred to as \textit{individual chance constraints} (e.\,g., see \cite{Li_etAl_2008_CC_PC_UncSys}).

\subsection{Problem Formulation} 

This paper considers  
the nonlinear system \eqref{eq:sys:f} with parametric uncertainties, and in which the initial states at sampling time $t_k$ are uncertain, e.\,g.\ due to uncertain state estimates.
Let $t_k$ denote the sampling time instances at which states $\bs{x}(t_k)$ become available.
Furthermore, denote by the vector $\bs{\pi}\in \mathbb{R}^{n_\pi}$ and function $\obs{u}$ a parametrization of the continuous-time input such that $\bs{u}(t) = \obs{u}(t,\bs{\pi})$, $t\in [ t_k, t_f]$.%
\footnote{For a piecewise-constant control input parameterization  $\bs{u}(t) = \bs{\pi}_i, \; t \in [t_{u,i}, t_{u,i+1}]$, $i=1,\ldots,N_u$, 
partition the time horizon $[t_k, t_f]$ into $N_u$ subintervals $[t_{u,i}, t_{u,i+1}]$ with $t_{k} = t_{u,1} < \ldots < t_{u,N_u} < t_{u,N_u+1} = t_f$.}

This work considers the following main problem.
\\\noindent\textbf{Finite-horizon stochastic nonlinear MPC with joint chance constraints:}
\textit{At each sampling time $t_k$ the following stochastic optimal control problem is solved
\begin{subequations}
\label{eq:prob1}
\begin{alignat}{5}
\underset{\bs{\pi}}{\text{min\ }} & J(\obs{x}(\cdot), \obs{u}(\cdot)) \label{eq:prob1:cost}\\
\text{subject to:\ } & \dot{\obs{x}}(t) = \bs{f}(\obs{x}(t),\obs{u}(t, \bs{\pi}),\bs{p}), &\forall t\in [t_k,t_f] \label{eq:prob1:dyn}\\ 
& \Prob{\bs{g}(\bs{X}) \le \bs{0} } \ge \beta, \label{eq:prob1:cc}\\
& \obs{u}(t,\bs{\pi}) \in \mathcal{U}, &\forall t \in  [t_k,t_f]\label{eq:prob1:U}\\ 
& \overline{x}_i(t_k) \sim \mu_{x_i(t_k)}, & i = 1,\ldots, n_x \label{eq:prob1:unc:x0}\\
& \bs{p}_i \sim \mu_{p_i}, & i = 1,\ldots, n_p \label{eq:prob1:unc:p}
\end{alignat}
\end{subequations}
where $\bs{\pi}$ denote the decision variables, $\obs{x}(t)$ denotes the states predicted by the nonlinear system model, and $\mu_{{x_i}(t_k)}$ denotes the probability distributions of the states at time $t_k$.
}

To facilitate closed-loop control, the stochastic optimal control problem~\eqref{eq:prob1} is often implemented in a receding-horizon mode or shrinking-horizon mode (cf.\ footnote \ref{fn:horizon}). 
The closed-loop control that is applied to the system~\eqref{eq:sys:f} is defined by the optimal solution $\bs{\pi}^\ast$ of \eqref{eq:prob1} at the sampling instants: $\bs{u}(t) = \obs{u}(t,\bs{\pi}^\ast)$, $t\in[t_k,t_{k+1}]$.


This work presents a framework to efficiently solve \eqref{eq:prob1}. 
In particular, the following problems are addressed.
\\\noindent\textbf{Problem 1: Propagation of the time-invariant probabilistic uncertainties \eqref{eq:prob1:unc:x0} and \eqref{eq:prob1:unc:p} through the nonlinear system dynamics \eqref{eq:prob1:dyn}.}\\
The problem is addressed using the polynomial chaos approach as presented in Sec.~\ref{sec:PC}.
This allows to efficiently sample from the probability distributions of the states and cost functions  (Sec.~\ref{sec:CC}) to address the next problem.
\\\noindent\textbf{Problem 2: Efficient evaluation of the chance constraints \eqref{eq:prob1:cc}.}\\
Accurate prediction of uncertain system dynamics typically requires a large number of samples, which can still be prohibitive even in the polynomial chaos approach. 
Moreover, sample-based approximations of \eqref{eq:prob1} may yield solutions that are infeasible for the original problem \eqref{eq:prob1} \cite{Calafiore_Campi_2005_MathProgr_UncConvProg_RA_ConfidenceLevels}. 
To reduce the risk of infeasibility due to the statistical error made due to the finite sample size while avoiding excessive sampling, the constraints can be tightened to make the entire problem more robust.
This is addressed in the following problem.
\\\noindent\textbf{Problem 3: Tightening the chance constraints \eqref{eq:prob1:cc} by $\beta_\text{cor} > \beta$ depending on the sample-size to guarantee a desired feasibility probability.}\\
A solution to this problem is presented in Sec.~\ref{sec:CC}.
The last problem addresses the efficiency of solving the stochastic optimal control problem \eqref{eq:prob1} using gradient-based optimization algorithms.
\\\noindent\textbf{Problem 4: Determining the gradients of the chance constraints \eqref{eq:prob1:cc}.}\\
A solution to the latter problem is presented in Sec.~\ref{sec:gradCC} and the overall framework is demonstrated in Sec.~\ref{sec:example}.

\section{Uncertainty Propagation for Nonlinear Systems Using Polynomial Chaos}
\label{sec:PC}

This work uses the polynomial chaos framework to solve Problem~1.
In the PC framework, spectral expansions in terms of orthogonal polynomials are used to represent stochastic variables and parameters \cite{Ghanem_Spanos_1991_SFE,wie38,xiu02,Kim_etal_2013_CSM__PCE_review}.
This allows deriving surrogate models, based on which the moments of the system states of the original system can be readily obtained. As shown in Sec.~\ref{sec:CC}, PC also allows for the sampling-based approximation and evaluation of chance constraints in a computationally efficient manner.
 
\subsection{Polynomial Chaos Expansion}
\label{sec:PCE}	

In the following, we assume that all uncertain parameters and uncertain initial conditions are functions of the standard random variables $\bs{\xi} \in \mathbb{R}^{n_\xi}$, which is denoted by $\bs{x}(t, \bs\xi)$, $\bs{p}(\bs\xi)$, etc.
The random variables $\xi_j$, $j = 1, \ldots, n_\xi$ are assumed to be independent with known probability distribution functions (PDFs) $\mu_{\xi_j}$, such that $\xi_j \in {L}^2(\Omega,\mathcal{F},\mu)$.
Here, ${L}^2(\Omega,\mathcal{F},\mu)$ is the Hilbert space of all random variables $\xi_j$ with finite variance $\Variance{\xi_j}<\infty$.
The triple $(\Omega,\mathcal{F},\mu)$ denotes the probability space on the basis of the sample space $\Omega$, $\sigma$-algebra $\mathcal{F}$ of subsets of $\Omega$, and probability measure $\mu$ on $\Omega$.

Le $v(t,\bs\xi)$ represent a state variable $x_i(t,\bs\xi)$ or any other (possibly) nonlinear function such as $g_i(\bs{x}(t,\bs\xi))$ in the chance constraints \eqref{eq:cons:CC}.
To explicitly derive the dependence of $v$ on the random variables $\bs\xi$, the following polynomial chaos expansion can be used
\cite{wie38,xiu02,Cameron_Martin_1947_Orthogonal_NonLin_Func}
\begin{align}
\label{eq:PCE}
v(t,\bs{\xi}) =
\sum_{\bs{\alpha}_i\in\mathcal{I}_\infty} \widetilde{v}_{\bs{\alpha_i}}(t)\Psi_{\bs{\alpha_i}}(\bs\xi).
\end{align}
The variables $\widetilde{v}_{\bs{\alpha_i}}$ denote the (deterministic) coefficients of the PC expansion (PCE), and $\Psi_{\bs{\alpha_i}}(\bs\xi)$ denote multivariate polynomials in the random variables $\bs\xi$ of total degree $\sum_{j=1}^{n_\xi} (\bs{\alpha_i})_j$.
The multivariate polynomials can be written as products of univariate polynomials: $\Psi_{\bs{\alpha_i}}(\bs\xi) \coloneqq \prod_{j=1}^{n_\xi}\Phi^{(\bs{\alpha_i})_j}_{\xi_j}$.
The polynomial $\Phi^m_{\xi_j}$ of the random variable $\xi_j$ is of degree $m$, where $\bigl\{\Phi^m_{\xi_j}\bigr\}_{m=0}^{P}$,  $j = 1,\ldots, n_\xi$ is an orthogonal basis\footnote{Such orthogonal bases are readily available for well-known standard distributions such as Normal, Uniform or Beta distributions \cite{xiu02}, or, in general, can be constructed for any distribution using moments \cite{Oladyshkin_Nowak_2012_Arbitrary_PCE} or Gram-Schmidt orthogonalization \cite{Gerritsma_etAl_2010_TD_PCE}.}
 with respect to the corresponding probability measures $\mu(\xi_j)$.

In \eqref{eq:PCE}, $\bs{\alpha_i}$ is the $i^\text{th}$ multi-index vector from the set
\begin{align*}
\mathcal{I}_P \coloneqq \left\{ \bs\alpha : \bs\alpha \in \mathbb{N}_{\ge 0}^{n_\xi}, \sum_{j=1}^{n_\xi} \alpha_j \le P \right\}
\end{align*}
with $P = \infty$.
For practical reasons, the infinite (weighted) sum of polynomials in Eq.~\eqref{eq:PCE} is truncated after $ \widetilde{P} \coloneqq \frac{(n_\xi+P)!}{n_\xi!P!}$ terms, where $P$ is called the order of the PC expansion.
The truncation can be written in a compact form as
\begin{align}
\label{eq:PCE:trunc}
v(t,\bs\xi) \approx \tbs{ v}(t)^\transp\bs{\Psi}(\bs{\xi}),
\end{align}
with 
\begin{align}
\label{eq:PCE:coeffs}
\tbs{ v}(t) \coloneqq \left[\widetilde v_{\bs{\alpha_1}}(t), \widetilde v_{\bs{\alpha_2}}(t),\ldots, \widetilde v_{\bs{\alpha_{\widetilde{P}}}}(t)\right]^\transp \in \mathbb{R}^{\widetilde{P}}
\end{align}
being the vector of coefficients for the PC expansion of variable $v$, and
\begin{align}
\label{eq:PCE:Phi_vec}
\bs{\Psi}(\bs\xi) \coloneqq \left[\Psi_{\bs{\alpha_1}}(\bs\xi),\Psi_{\bs{\alpha_2}}(\bs\xi),\ldots, \Psi_{\bs{\alpha_{\widetilde{P}}}}(\bs\xi)\right]^\transp \in \mathbb{R}^{\widetilde{P}}
\end{align}
being the vector of the multivariate polynomials.

The next step is to determine the values of the PC expansion coefficients \eqref{eq:PCE:coeffs}, which will be needed to approximate the probability distributions for the evaluation of the chance constraints (see Sec.~\ref{sec:CC}).
Two different approaches to determine the PC expansion coefficients are explained in the next two subsections.


\subsection{Determination of the PC Expansion Coefficients Using Collocation}
\label{sec:PCE:collocation}

The PC expansion coefficients can be determined using the so-called probabilistic collocation methods (e.g., see \cite{Nagy_Braatz_2007_JPC_UncertAnalysis_PCE,fag12,Mesbah_etAl_2014_IFACWC__aFDI_ProbUnc,Tatang_etAl_1997_Efficient_ParamUnc_Collocation} and references within). In the collocation methods, $N_{MC}$ samples are drawn from the known distributions of uncertainties and, subsequently, are used to solve the nonlinear process model~\eqref{eq:sys:f}. The PC expansion coefficients can then be obtained in a least squares sense through minimizing the residuals between the PC expansion and the nonlinear model predictions $v(t,\bs{\xi^{[j]}})$. Note that an explicit analytic solution to the resulting unconstrained optimization problem ${{\displaystyle\min} \atop {\widetilde{\bs v}(t)}} \sum_{j=1}^{N_{MC}} \left( \widetilde{\bs v}(t)^\transp \bs{\Psi}(\bs{\xi^{[j]}})  - v(t,\bs{\xi^{[j]}}) \right)^2$ is readily available.

The complexity of process dynamics in terms of nonlinearities may require a large number of samples $N_{MC}$ to obtain adequate estimations of the PC expansion coefficients. This may render real-time control applications computationally prohibitive. Next, an alternative approach is presented for determining the PC expansion coefficients that relies on the solution of an extended set of ordinary differential equations.  

\subsection{Determination of the PC Expansion Coefficients Using Galerkin Projection}
\label{sec:PCE:Galerkin}

In the following it is assumed that the considered system is polynomial in the states, i.\,e. the dynamics for the $i^\text{th}$ state
\begin{align}
\label{eq:sys:polynomial}
\dot{x}_i(t, \bs\xi) = f_i(\bs x(t, \bs\xi),\bs u(t), \bs p(\bs\xi)), \quad\forall i = 1, \ldots, n_x
\end{align}
is of the form $c \prod_{l=1}^{n_x}x_l^{\gamma_{l}}(t, \bs\xi) \sigma_u(\bs u(t))\sigma_p(\bs p(\bs\xi))$, where $c$ is a constant and $\gamma_l$ is the degree of variable $x_l$.
$\sigma_p$ and $\sigma_u$ are nonlinear functions of the parameters and inputs, respectively.
Note that such a representation can always be obtained exactly by state-lifting or immersion \cite{Ohtsuka_2005_IEEETAC_Immersion} under the conditions given in Sec.~\ref{sec:problem}.
Such a polynomial model structure enables explicit derivation of differential equations for the PC coefficients $\tbs{x}_i$ using Galerkin projection \cite{Ghanem_Spanos_1991_SFE}.
The Galerkin projection requires evaluation of multi-dimensional integrals, which can be solved exactly and efficiently offline for polynomial systems in the form \eqref{eq:sys:polynomial}.
For generality, the (uncertain) parameters need not appear polynomially in \eqref{eq:sys:polynomial}.
If $\sigma_p$ is not polynomial, a PC approximation of $\sigma_p$ can be determined using the collocation approach described in the previous section \ref{sec:PCE:collocation}.

\subsubsection{Galerkin Projection of the System Dynamics}
\label{sec:PCE:Galerkin:dyn}

To obtain the PC coefficients for the states $\tbs{x}_i(t) $, an extended system of ordinary differential equations is derived using Galerkin projection \cite{wie38,xiu02}.
This system is obtained by formal derivative of the PC expansion \eqref{eq:PCE:trunc} for the state variables (yielding ${\dot{\tbs{ x}}_i(t)}^\transp\bs{\Psi}(\bs{\xi})$) and by inserting the PC expansion of the state variables \eqref{eq:PCE:trunc} into the system dynamics \eqref{eq:sys:polynomial}.
Subsequently, the projection is performed by computing the inner products $\int f_i(\tbs{x},\bs{u}(t),\tbs{p}) \Psi_{\bs{\alpha_i}}(\bs\xi) \mu_1(d\xi_1)\cdots\mu_{n_\xi}(d\xi_{n_\xi})$ of the resulting equation and the different polynomials in \eqref{eq:PCE:Phi_vec}.
By employing orthogonality, this results in the set of ordinary differential equations (ODEs) describing the dynamics of the coefficients
\begin{equation}
\label{eq:sys:PCE}
\dot{\tbs{ x}}_i(t) = 
\tbs{ f}_i(\tbs{ x}_1(t),\ldots,\tbs{ x}_{n_x}(t) ,\bs u(t),\tbs{p}_1,\ldots,\tbs{p}_{n_p}), \\ \forall i = 1, \ldots, n_x,
\end{equation}
where $\tbs{x}_{i}$ and $\tbs{p}_i$ are the vectors of coefficients of the PC expansions (\cf \eqref{eq:PCE:coeffs}) of the states and parameters, respectively. 
The system \eqref{eq:sys:PCE} has extended state space dimension $\widetilde{P}\,n_x$ and describes the dynamics of the PC expansion coefficients.
Thus, by projection onto the orthogonal polynomials, 
the explicit dependencies on the random variables are eliminated.
The solution to this set of ODEs  can then be used for efficient sample evaluations (\cf Sec.~\ref{sec:PCE:approx}).

To compute the inner products, multi-dimensional integration is required.
Even though this is in general cumbersome, it is particularly easy for polynomial systems \eqref{eq:sys:polynomial}.
In this case, the integral can be efficiently and exactly computed using Gauss Quadrature \cite{Gautschi_etAl_2004_Book__Orthogonal_Polynomials}.
In addition, due to (power) orthogonality \cite{Gautschi_etAl_2004_Book__Orthogonal_Polynomials} most projection integrals ($\sim\!97\%$) are zero such that the computational burden can be reduced significantly.
For more details and the structure of \eqref{eq:sys:PCE}, see \cite{Streif_etAl_2014_IFACWC__PCE_ExpDesign_MD_PC}.

\subsubsection{Galerkin Projection of the Initial Conditions}

Once \eqref{eq:sys:PCE} has been determined as described in Secs.~\ref{sec:PCE:collocation} and \ref{sec:PCE:Galerkin}, the initial conditions $\tbs{x}_i(t_k)$ are needed for numerical solution of the set of ODEs. Since the initial conditions $\bs{x}(t_k)$ are assumed to be functions of the uncertainties $\bs\xi$, they can be obtained by projection of the corresponding PC expansion \eqref{eq:PCE} onto the different orthogonal polynomials $\widetilde{\Psi}_{\bs{\alpha_i}}(\bs\xi)$
\begin{equation}
\label{eq:PCE:IC}
\widetilde{x}_{i,j}(t_k) = 
\frac
{\int \bs{x}(t_k,\bs\xi) \Psi_{\bs{\alpha_i}}(\bs\xi)\mu_1(d\xi_1)\cdots\mu_{n_\xi}(d\xi_{n_\xi})}
{\int \left(\Psi_{\bs{\alpha_i}}(\bs\xi)\right)^2\mu_1(d\xi_1)\cdots\mu_{n_\xi}(d\xi_{n_\xi})},\\ 
\forall i = 1, \ldots, n_x, \forall j = 1, \ldots, \widetilde{P}.
\end{equation}

\subsection{Computation of the Moments and Efficient Sampling to Approximate Probability Distributions}
\label{sec:PCE:approx}

Once the PC expansion coefficients \eqref{eq:PCE:coeffs} are determined, the probability density of $v(t,\bs{\xi})$ can be approximated using sampling.
With that, the approximation of the probability distributions of stochastic variables $v$ or the evaluation of chance constraints can be done efficiently as shown in Sec.~\ref{sec:CC}.

Let $\bs{\xi}^{[j]}$ be samples drawn from the multivariate distributions of $\bs\xi$.
Then $v(t,\bs{\xi^{[j]}_i})$ is obtained from \eqref{eq:PCE:trunc} by evaluating the orthogonal polynomials \eqref{eq:PCE:Phi_vec} for the different samples $\bs{\xi^{[j]}}$, and by vector multiplications with the vector of the PC expansion coefficients obtained from the solution of \eqref{eq:sys:PCE}.

Besides such a sampling-based approach, the probability distributions can be approximated based on moments (see discussion and references in \cite{Streif_etAl_2014_IFACWC__PCE_ExpDesign_MD_PC}).
The moments can be determined directly from polynomial chaos expansions without further approximations (see \cite{Fisher_Bhattacharya_2009_LQR_PCE,Streif_etAl_2014_IFACWC__PCE_ExpDesign_MD_PC}), which is computationally cheap for low-order moments.
The moments are also required to compute the cost function \eqref{eq:cost} in \eqref{eq:prob1}.

\section{Sample Approximation of Chance Constraints with Guaranteed Feasibility Probability}
\label{sec:CC}

In this section, Problems~2 and 3 are addressed.
In particular, the satisfaction probability of the joint chance constraint \eqref{eq:cons:CC} is approximated using samples generated with the PC framework. 
Even though a large number of samples can be generated in a computationally efficient manner using the PC expansion (see Sec.~\ref{sec:PC}), the evaluation of the nonlinear functions $g_i\left(\tbs{X}, \bs{\xi}\right)$, $i =1,\ldots, N_g$ and their gradients (see Sec.~\ref{sec:gradCC}) may still be prohibitive for a large sample size or large $N_g$.
To increase the computational efficiency, it is therefore desired to evaluate as few samples as possible.
However, fewer samples increase the width of the confidence interval (i.\,e.\ reduce the quality) of the estimate of the satisfaction probability.
A low confidence bears the risk that a feasible solution to the sample-approximated chance constrained optimization problem is infeasible for the original problem (e.\,g.\ for a different or larger set of samples) \cite{Calafiore_Campi_2005_MathProgr_UncConvProg_RA_ConfidenceLevels}.
To increase the \textit{confidence level of feasibility}, which we call \textit{feasibility probability}, one can make the controller more robust by tightening the chance constraints (i.\,e., using a $\beta_\text{cor} > \beta$ in Eq.~\eqref{eq:cons:CC}).
This section proposes a statistical analysis to determine a-priori the 
the constraint tightening $\beta_\text{cor}$ for which a desired feasibility probability $1-\alpha$ ($\alpha \in (0,1)$) and a desired robustness level (satisfaction probability $\beta$) can be guaranteed.

Note that this section considers the PC expanded system \eqref{eq:sys:PCE} rather than the original nonlinear system \eqref{eq:prob1:dyn}.

\subsection{Satisfaction Probability}

Suppose that the inequalities in the chance constraints \eqref{eq:cons:CC} are expressed using the PC expansion
\begin{align}
\label{eq:cons:PCE:CC}
\Prob{\bs{g}\left(\tbs{X}, \bs{\xi}\right) \le \bs{0} } \ge \beta,
\end{align}
with $\beta$ as in \eqref{eq:cons:CC}, and 
\begin{equation*}
\tbs{X} \coloneqq \Bigl[ \tbs{x}_1(t_{x,1})^\transp, \ldots, \tbs{x}_{n_x}(t_{x,1})^\transp,\ldots,\\ \tbs{x}_1(t_{x,N_x})^\transp,\ldots,
\tbs{x}_{n_x}(t_{x,N_x})^\transp\Bigr]^\transp
\end{equation*}
are given from the simulation of the PC expanded system \eqref{eq:sys:PCE}.
The probability of satisfaction in \eqref{eq:cons:PCE:CC} is given by
\begin{equation}
\label{eq:prob}
\Prob{\bs{g}\left(\tbs{X}, \bs{\xi}\right) \le \bs{0} }
\coloneqq \\\lint_{-\infty}^{+\infty}\lint_{-\infty}^{+\infty}\cdots\lint_{-\infty}^{+\infty}
I_{\mathcal{G}}\left(\bs{\xi}\right)\, \mu_{1}(d\xi_1)\mu_{2}(d\xi_{2})\cdots\mu_{n_\xi}(d\xi_{n_\xi}).
\end{equation}
$I_{\mathcal{G}}$ is the indicator function
\begin{align*}
I_{\mathcal{G}}\left(\bs{\xi}\right)
\coloneqq
\begin{cases}
1 & \text{if\ } \bs{\xi} \in \mathcal{G}\\
0 & \text{otherwise},
\end{cases}
\end{align*}
where $\mathcal{G} \in \mathbb{R}^{n_\xi}$ denotes the set where all constraints are satisfied
\begin{align*}
\mathcal{G} \coloneqq \left\{ \bs{\xi}\in \mathbb{R}^{n_\xi} \ |\ g_i\left(\tbs{X},\bs{\xi}\right)\le 0, \forall i = 1,\ldots,N_g \right\}.
\end{align*}

Evaluation of the multidimensional integral in \eqref{eq:prob} is in general very difficult due to the non-convexity of the set $\mathcal{G}$ and the weighting with respect to the probability measures $\mu(\xi_i)$, $i=1,\ldots,n_\xi$.

In this work, sampling is used to approximate the probability of satisfaction \eqref{eq:prob}.
For this purpose, $n_\xi$-dimensional independent and identically distributed samples $\bs{\xi}^{[1]}, \ldots, \bs{\xi}^{[N_S]}$ are drawn from the distributions $\mu_i(\xi_i)$, $i=1,\ldots,n_\xi$.
The sample-average approximation of $\Prob{\bs{g}\left(\tbs{X}, \bs{\xi}\right) \le \bs{0} }$ is given by
\begin{align}
\label{eq:prob:approx}
\Prob{\bs{g}\left(\tbs{X}, \bs{\xi}\right) \le \bs{0}}
\approx
\frac{1}{N_S}
\sum_{j=1}^{N_S}I_{\mathcal{G}}\left(\bs{\xi}^{[j]}\right).
\end{align}

\subsection{Required Sample Size and Constraint Tightening for a Guaranteed Feasibility Probability}
\label{sec:CC:confidence}

The main results of this section are the following theorems that allow determining systematically the constraint tightening and sample size for which desired confidence level in the approximation \eqref{eq:prob:approx} of \eqref{eq:prob} is guaranteed. 
\\[1ex]\noindent\textbf{Theorem~1 (Constraint Tightening and Confidence in the Sample Approximation of the Chance Constraints):}
\textit{
If $\beta_\text{cor} > \beta$ is chosen such that 
\begin{align}
\label{eq:N_S}
1- \text{betainv}\left(1 - \frac{\alpha}{2},N_S +1 - \floor{\beta_\text{cor} N_S}, \floor{\beta_\text{cor} N_S}\right) \ge \beta,
\end{align}
then the sample approximation \eqref{eq:prob:approx} of the tightened chance constraints 
$\Prob{\bs{g}\left(\tbs{X}, \bs{\xi}\right) \le \bs{0} } \ge \beta_{\text{cor}}$ using $N_S$ samples guarantees a confidence level of $1-\alpha$ of 
the original chance constraint $\Prob{\bs{g}\left(\tbs{X}, \bs{\xi}\right) \le \bs{0} } \ge \beta$ (\eqref{eq:cons:PCE:CC} resp.\ \eqref{eq:prob1:cc}).
Furthermore, $\beta_\text{cor}$ is the least conservative constraint tightening that can be chosen.
\\The operator $\floor{\cdot}$ denotes rounding towards $-\infty$, and $\text{betainv}$ denotes the inverse of the cumulative Beta-distribution.
}\\
\noindent\textbf{Proof:} \textit{The proof uses standard results in statistics and is presented in Appendix~A.\hfill$\square$}\\[1ex]
\indent%
The theorem tightens the chance constraints $\beta_\text{cor}$ to compensate for the statistical error made due to the finite sample size $N_S$.
To this end, it employs the lower bounds of the confidence interval obtained from a statistical analysis (details see the proof).
Note that the analysis is independent of the specifics of the optimization problems, in particular of the chance constraints and the system dynamics.
Note also that the analysis neither depends on the number of decision variables nor requires convexity assumptions as in \cite{Calafiore_Campi_2006_IEEETAC__Scenario_Approach_RobustControl,Campi_Garatti_2010_JOTA__SA_CC,Calafiore_2010_SIAMJO__Random_ConvexPrograms,
Garatti_Campi_2013_CSM__Modulating_Robustness_ControlDesign,
Campi_Garatti_2008_SIAMJO__ExactFeasibility_RA_ConvexPrograms}.
Furthermore, bounds on the sample complexity have been presented in \cite{Alamo_etAl_2010_Sample_Complexity_Prob_Analysis_Design,Alamo_etAl_2014_RA_UncSys_Sample_Complexity} in a different context.

\begin{figure}[th]
\centering
\input{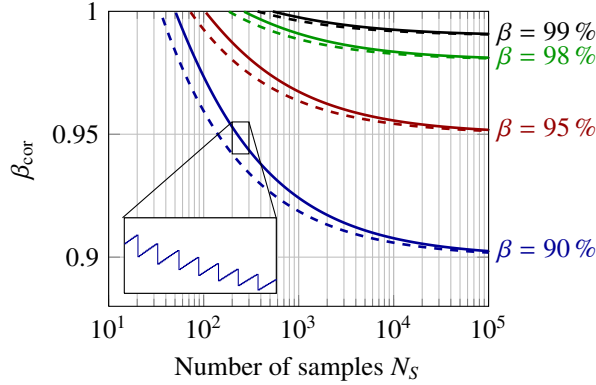}
\caption{Sample size $N_S$ and corrected satisfaction probability $\beta_{\text{cor}}$ (\ie chance constraint tightening) for different required satisfaction probabilities $\beta$ for a feasibility probability of $1-\alpha = 0.99$ (solid lines) and $1-\alpha = 0.95$ (dashed lines). 
The inset shows the non-monotonicity of the curves.}
\label{fig:CC:minSamples}
\end{figure}

Due to the independence of the analysis on the specifics of the optimization problems and constraints, Theorem~1 can be applied offline (and needs to be done only once) to systematically satisfy prespecified probabilistic requirements.
The implicit relationships of \eqref{eq:N_S} are plotted in Fig.~\ref{fig:CC:minSamples} for selected values of $\alpha$, $\beta$, and $N_S$ as a reference, and can be derived for other values similarly.
It is noteworthy that more samples does not necessarily mean a tighter confidence interval \cite{Rossi_CC_ConfidenceInterval}, which can be also seen at the lack of monotonicity of the curves in the inset of Fig.~\ref{fig:CC:minSamples}.

Theorem~1 can now be used to guarantee a certain confidence level of the feasibility of a solution to \eqref{eq:prob1}.
\\[1ex]\noindent\textbf{Theorem~2 (Feasibility Probability):}
\textit{Consider the chance constrained stochastic optimal control problem \eqref{eq:prob1} and let $\bs{\pi}'$ be a feasible solution to the sample-approximated form of problem \eqref{eq:prob1} for a given $\beta_\text{cor}$ according to Theorem~1.
Then $\bs{\pi}'$ is a feasible point of \eqref{eq:prob1} with probability not less than $1-\alpha$.
}\\
\noindent\textit{\textbf{Proof:} See Appendix~B.\hfill$\square$}\\[1ex]

Theorem~2 guarantees that a solution found by sample approximation is also a solution to the original problem \eqref{eq:prob1} (i.\,e.\ for an infinite number of samples) with the specified confidence level $1-\alpha$.

\section{Sample Approximation of the Gradients for Efficient Optimization}
\label{sec:gradCC}

\begin{figure*}[tp!]
{\centering
\includegraphics[width=0.48\linewidth]{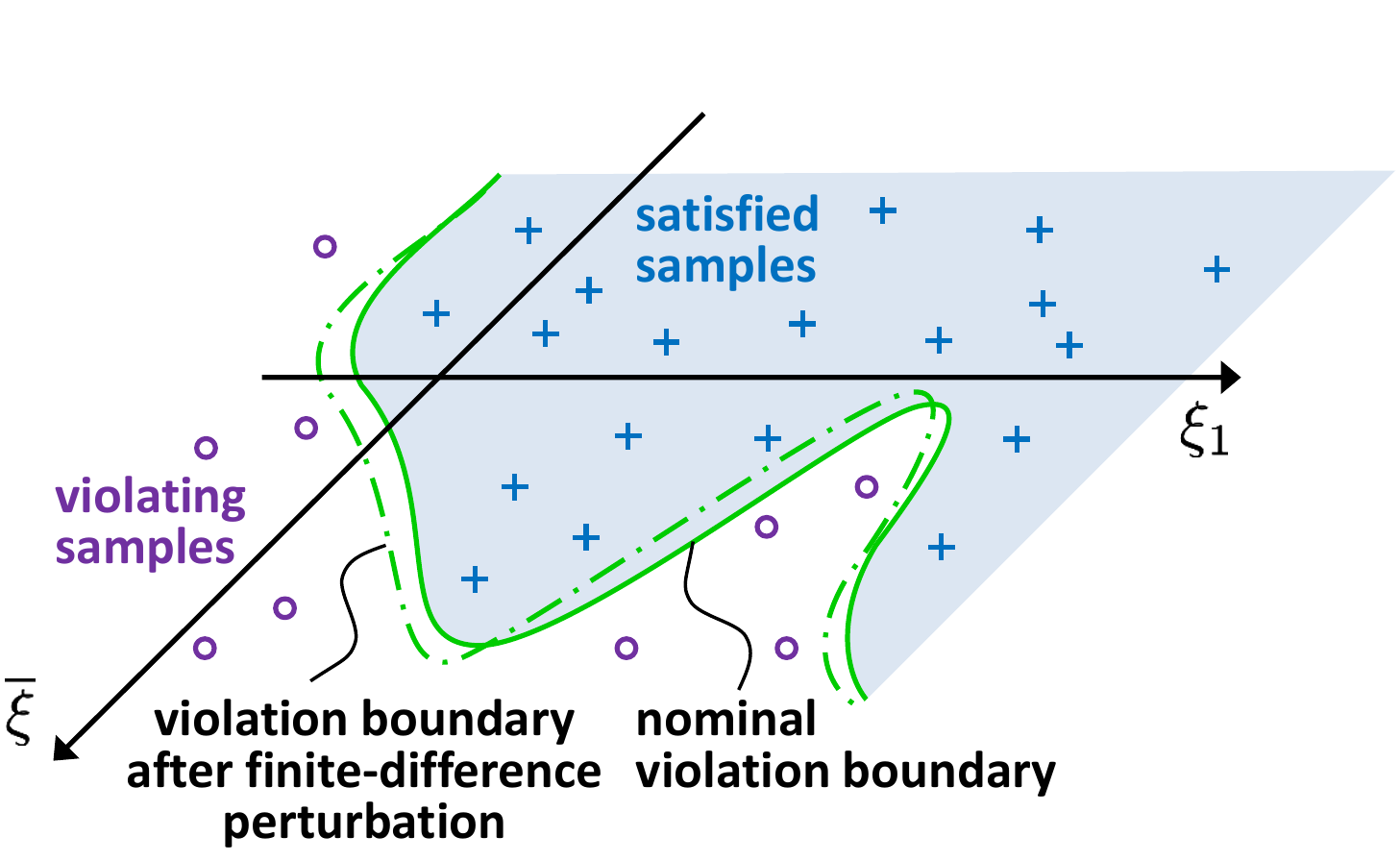}
\includegraphics[width=0.48\linewidth]{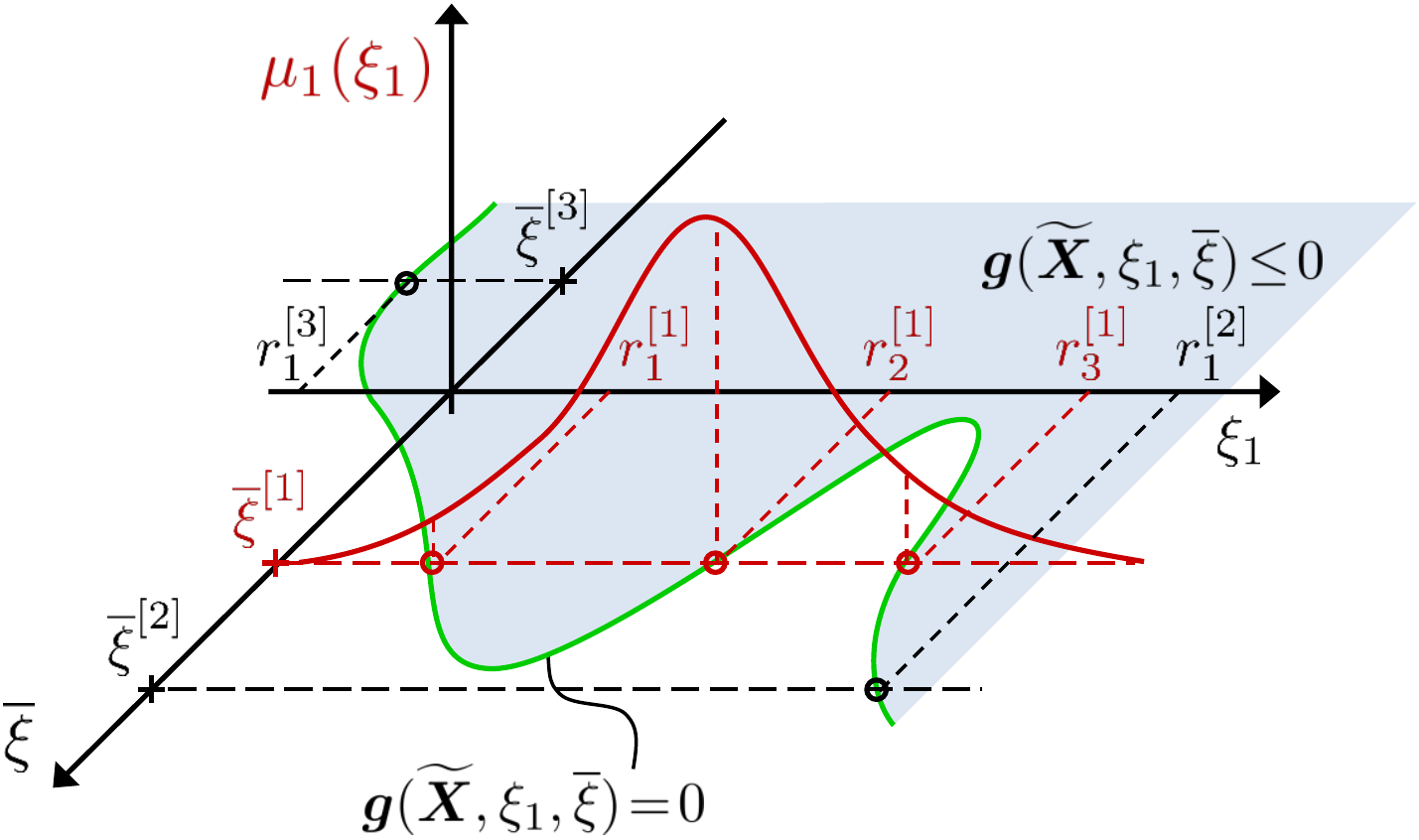}\\}
(a)\hspace{0.48\linewidth}(b)
\caption{(a) Samples satisfying ${g}(\cdot)\le0$ are shown by blue $\bs{+}$, and violating samples by purple $\circ$.
The satisfaction probability is given by ratio ``number of satisfying samples''/``total number of samples''. This ratio changes discontinuously due to finite-difference perturbations $\bs{\widetilde{X}}+\bs{\Delta \widetilde{X}}$ because a discrete number of samples change their validity. Note that  the ratio may not even change at all if the validity of none of the samples is affected.
(b) The samples $\bar{\xi}^{[1]}, \bar{\xi}^{[2]}, \bar{\xi}^{[3]}$ (symbol $\bs{+}$) are drawn at random from the distributions $\mu_2(\xi_2)$, where $\bar{\xi}=\xi_2$.
The equation ${g}(\tbs{X},\xi_1,\bar{\xi}^{[j]})={0}$ defines the violation boundary (shown in green) and it is solved for $\xi_1$.
For $\bar{\xi}^{[1]}$, this yields three solutions $\bs{r}^{[1]} = [r_1^{[1]}, r_2^{[1]},r_3^{[1]}]$ (shown in red).
A different number of (real-valued) solutions is obtained for the other samples as shown for the red and black dashed lines.
Piecewise integration of the probability measure $\mu_1(\xi_1)$  is performed with integration limits given by the elements of $\bs{r}^{[1]}$.
Note that the integral between the limits $r_1^{[1]}$ and $r_2^{[1]}$ contributes to the approximation of the gradient only if the constraint  $g(\tbs{X},\hat{r}_1^{[1]},\bar\xi)\le0$ for an arbitrary point $\hat{r}_1^{[1]}$  (not shown) between the integration limits.}
\label{fig:CC:approach}
\end{figure*}

This section proposes a solution to Problem~4, that is a sample-based approximation of the analytic expressions of the gradients of the chance constraints.
The proposed approach avoids approximations of the gradients by finite difference methods, which typically slows down gradient-based optimization algorithms (see \eg \cite{Imsland_etAl_2010_MB_Contr_Comparison_FiniteDiff_SensODE}).
Furthermore, finite-differencing of the chance constraints can leads to poor estimates due to the discrete nature of sample-approximations as explained in Fig.~\ref{fig:CC:approach}a and \cite{Garnier_etAl_2009_Approximative_Gradients_CC,Royset_Polak_2004_SAA_CC}.
\footnote{Note that more sophisticated sampling methods such as importance sampling could be used at the cost of much higher computational demands.}

The gradients of the joint chance constraints with respect to the input parametrization $\bs{\pi}$ are formally given by
\begin{align}
\label{eq:gradCC:parts}
\frac{d\Prob{\bs{g}\left(\tbs{X}, \bs{\xi}\right)\le \bs{0}}}{d\bs{\pi}} 
= \frac{\partial \Prob{\bs{g}\left(\tbs{X}, \bs{\xi}\right)\le \bs{0}}}{\partial \tbs{X}}
\frac{\partial \tbs{X}}{\partial \bs{\pi}}.
\end{align}

Note that the first-order sensitivities $\frac{\partial\tbs{X}}{\partial\bs{\pi}}$ are obtained from the solution of the sensitivity equations, which are almost always needed to speed up gradient-based optimization algorithms. 

The main result of this section is the following proposition, namely sample-based approximation of the analytic gradients of \eqref{eq:prob} with respect to the input parametrization variables $\bs{\pi}$.
Note that the dependence of the states $\tbs{X}$ on $\bs{\pi}$ is not explicitly written to shorten the notation.
The derivation is illustrated in Fig.~\ref{fig:CC:approach}b and the constructive proof is given with all technical details in the Appendix~C. 
\\[1ex]\noindent\textbf{Proposition~1 (Gradients of the Sample-Approximated Chance Constraints):}
\textit{%
Suppose the following to be given:
\begin{enumerate}
\item first-order sensitivities $\frac{\partial\tbs{X}}{\partial\bs{\pi}}$;
\item $\obs{\xi}^{[j]}$, $j=1, \ldots, N_S$ independent and identically distributed samples from the vector $\obs{\xi}\coloneqq [\xi_2, \ldots, \xi_{n_\xi} ]$;
\item for each sample  $\obs{\xi}^{[j]}$ the vector $\bs{r}^{[j]} \coloneqq [r_0^{[j]}, r_1^{[j]}, \ldots, r_{n_r^{[j]}+1}^{[j]}]$ with unique elements and sorted in ascending order; the vector is defined by the solutions of $\bs{g}(\tbs{X},\xi_1,\bar{\xi}^{[j]})=0$ with respect to $\xi_1$; $r_0^{[j]}$ ($r^{[j]}_{n_r^{[j]}+1}$) is the minimum (resp.\ maximum) of the support of $\xi_1$.
\end{enumerate}
Then a sample-based approximation of $\frac{\partial \Prob{\bs{g}\left(\tbs{X}, \bs{\xi}\right)\le \bs{0}}}{\partial \tbs{X}}$ in \eqref{eq:gradCC:parts} is given by
\begin{equation}
\label{eq:gradProb:h}
\frac{\partial \Prob{\bs{g}\left(\tbs{X}, \bs{\xi}\right) \le  \bs{0} }}
{\partial \tbs{X}}
\approx\\
\frac{1}{N_S}\sum_{j=1}^{N_S}
\Bigg[
\sum_{i=0}^{n_r^{[j]}}
I_{\mathcal{G}}\left(\frac{r^{[j]}_i + r^{[j]}_{i+1}}{2},\obs{\xi}^{[j]}\right)
\\
\left(
\mu_1\left(r_{i+1}^{[j]}\right)
\pfrac{r_{i+1}^{[j]}}{\tbs{X}}
-\mu_1\left(r_{i}^{[j]}\right)
\pfrac{r_{i}^{[j]}}{\tbs{X}}
\right)
\Bigg],
\end{equation}
in which
\begin{align}
\label{eq:gradProb:h:implicitFunc}
\pfrac{r_{i}^{[j]}}{\tbs{X}} =
\begin{cases}
0, & i \in \{0, n_r^{[j]}+1\}\\
 -\pfrac{g\left(\tbs{X}, r_{i}^{[j]},\obs{\xi}^{[j]}\right)}{\tbs{X}} \left( \pfrac{g\left(\tbs{X}, r_{i}^{[j]},\obs{\xi}^{[j]}\right)}{\xi_1} \right)^{-1},
& i \in \{ 1,\ldots,n_r^{[j]}\}.
\end{cases}
\end{align}
}\\
\noindent\textbf{Proof:} \textit{See Appendix~C.\hfill$\square$}\\[1ex]

The vector $\bs{r}^{[j]}$ defines the limits of the integration of the PDF of the random variable $\xi_1$ (for further explanations see Fig.~\ref{fig:CC:approach} and the proof).
Note that it is required to solve polynomial equations $\bs{g}$ in $\xi_1$ to obtain the vector of integration limits. However, the polynomials can be solved either analytically for low order polynomials, or numerically using efficient root finding algorithms.

It is important to note that Proposition 1 makes the implicit assumption that only samples with unique roots $\bs{r}^{[j]}$ are used.
It is expected that samples with nonunique roots are relatively rare, such that the error will be negligible when such samples are discarded (in the case study in Sec.~\ref{sec:example} no samples were discarded).

With Proposition~1, the gradients can be computed efficiently since time-consuming finite differencing is avoided. Furthermore, the approach does not suffer from the discretization effects shown in Fig.~\ref{fig:CC:approach}a.
Note that the presented approach is inspired by \cite{Royset_Polak_2004_SAA_CC}.
However, the extension made in Proposition 1 allows for a much broader applicability such as polynomial and joint chance constraints, as well as non-gaussian probability distributions.

\section{Stochastic NMPC of a Williams-Otto Reactor}
\label{sec:example}

In this section, the solution to \eqref{eq:prob1} will be illustrated based on shrinking horizon SNMPC of a William-Otto semi-batch reactor.
To this end, uncertainties are propagated through the nonlinear system dynamics using the polynomial chaos approach presented in Sec.~\ref{sec:PC}. 
The sample-average approximation of the chance constraints \eqref{eq:prob:approx} and of their gradients \eqref{eq:gradProb:h} proposed in Sec.~\ref{sec:CC} resp.\ Sec.~\ref{sec:gradCC} are used.
The theorems presented in Sec.~\ref{sec:CC:confidence} are used to guarantee a desired feasibility probability by choosing the required constraint tightening and sample size according to equation \eqref{eq:N_S}.


\subsection{Williams-Otto Reactor}
\label{sec:example:WOR}

The Williams-Otto semi-batch reactor is considered,
in which the reactions $A + B \overset{\rho_1}{\rightarrow} C$, $C + B \overset{\rho_2}{\rightarrow} P + E$, and $P + C \overset{\rho_3}{\rightarrow} G$ take place with the associated reaction rates $\rho_1$, $\rho_2$, and $\rho_3$ \cite{Hannemann_Marquardt_2010_Adjoints_Hessians_Lagragian_DynOpt_WilliamsOttoReactor,Williams_Otto_1960_Reactor}.
The reactant $A$ is introduced into the reactor at the beginning of the batch process, whereas reactant $B$ is fed into the reactor with feed rate $u_1(t)$.
During the exothermic reactions, the products $P$ and $E$ as well as the side-product $G$ are formed.
The reactor temperature is directly manipulated using input $u_2(t)$.

The dynamic model consists of seven differential equations
\begin{equation}
\label{eq:WOR}
\arraycolsep=1pt
\begin{array}{rcrrrrr}
	\dot{x}_1(t) &=& - \rho_1(t) & & &&- \frac{x_1(t) u_1(t)}{ x_7(t)} \\
	\dot{x}_2(t) &=& - \rho_1(t) &- \rho_2(t) & &  + \frac{c_{B,in} u_1(t)}{ x_7(t)} &- \frac{x_2(t) u_1(t)}{ x_7(t)} \\
	\dot{x}_3(t) &=& \rho_1(t) &- \rho_2(t) &- \rho_3(t) &&- \frac{x_3(t) u_1(t)}{x_7(t)} \\	
	\dot{x}_4(t) &=& & \rho_2(t) & - \rho_3(t)  &&- \frac{x_4(t) u_1(t)}{ x_7(t)} \\
	\dot{x}_5(t) &=& & \rho_2(t)&  &&- \frac{x_5(t) u_1(t)}{ x_7(t)} \\
	\dot{x}_6(t) &=& & &\rho_3(t) &&-\frac{x_6(t) u_1(t)}{ x_7(t)} \\
	\dot{x}_7(t) &=& u_1(t),
\end{array}
\end{equation}
where $\bs{x}(t) \coloneqq [x_1(t), x_2(t), x_3(t), x_4(t), x_5(t), x_6(t)]^\transp$ denotes the molar concentrations (in units of $\text{mol}/\text{m}^3$) of $A$, $B$, $C$, $P$, $E$, and $G$ with respect to the volume $x_7(t)$ (in units of $\text{m}^3$) contained in the reactor at time $t$ (in seconds). $c_{B,in} = 5$ is the molar concentration of $B$ in the inlet reactor feed $u_1$ (in units of $\text{m}^3/\text{s}$) and $u_2$ (in units of $K$).
The last terms in the first six equations are the dilution effect due to inflow of reactant $B$.
The nonlinear state and input dependent reaction rates are
\begin{align*}
\rho_1(t) &= k_1 x_1(t) x_2(t)\exp(-6666.7/u_2(t))  \\
\rho_2(t) &= k_2 x_2(t) x_3(t)\exp(-8333.3/u_2(t))  \\
\rho_3(t) &= k_3 x_3(t) x_4(t)\exp(-11111/u_2(t)) .
\end{align*}
The parameter values of the reaction kinetics are taken from \cite{Hannemann_Marquardt_2010_Adjoints_Hessians_Lagragian_DynOpt_WilliamsOttoReactor,Williams_Otto_1960_Reactor} and have been converted to SI units. The initial conditions at the beginning of the batch are: 
\begin{equation*}
	\bs{x}(0)^\transp = \left[10\,\,0\,\,0\,\,0\,\,0\,\,0\,\,2 \right]. 
\end{equation*}


\subsubsection{Stochastic Optimal Control Problem}

The reaction rate constants $k_1$, $k_2$, and $k_3$ (in units of $\text{m}^3/(\text{mol}\,\text{s})$) are uncertain and probabilistically distributed according to Normal distributions.
\begin{align*}
k_1 &\sim \operatorname{Norm}(1.6599 \cdot 10^6, 1.6599 \cdot 10^5)\\
k_2 &\sim \operatorname{Norm}(7.2117 \cdot 10^8, 7.2117 \cdot 10^7)\\
k_3 &\sim \operatorname{Norm}(2.6745 \cdot 10^{12}, 2.6745 \cdot 10^{11}),
\end{align*}
where the first argument in $\operatorname{Norm}(\cdot,\cdot)$ specifies the mean and the second the variance. 
The variances are chosen to be $10\,\%$ of the mean values taken from \cite{Williams_Otto_1960_Reactor}. 

The objective of the process is to maximize the profit at the end of the batch, which is the difference
between the sales of the products $E$ and $P$ (1.0 monetary units per mole given by $x_4(t_f)x_7(t_f)$ resp.\ $x_5(t_f)x_7(t_f)$ with $x_7(t)$ being the volume at time $t$) and the costs of raw material $B$ ($0.5$ monetary units per mole).
The objective is defined by
\begin{multline*}
J(\bs{x}(\cdot),\bs{u}(\cdot)) = - 0.5 c_{B,in}\left( \Expect{x_7(t_f)} - \Expect{x_7(t_0)} \right) \\+  \Expect{x_5(t_f)x_7(t_f)}
+ 2\Expect{x_4(t_f)x_7(t_f)} \\- 10\Big(\Variance{x_5(t_f)x_7(t_f)}+\Variance{x_4(t_f)x_7(t_f)}\Big),
\end{multline*}
where the expectations and variances are computed as shown in Sec.~\ref{sec:PCE:approx}.
The objective takes into account the mean values of the amounts of the desired end-products $P$ and $E$, as well as the variances of the end products weighted by a positive constant (in units of $1/\text{mol}$) to reduce the variance of the end-products.

During the batch, constraints on the inputs must be satisfied at the time-points $t_{u,k}$ 
\begin{subequations}
\begin{align}
0 \le u_1(t) & \le 0.002 \\
313 \leq u_2(t) & \leq 363 \\
\left|u_2(t)-u_2(t)\right| & \le  1.
\end{align}
\end{subequations}
In addition, a constraint on the side-product $G$ ($x_6$) at the final time $t_f$ is defined to keep the amount of the undesired side-product at the end of the batch below a certain threshold.
This is because the batch products would be worthless (due to expensive filtering or purifications) if the threshold is exceeded.
Due to physical limitations, the volume $x_7$ at the final time $t_f$ is also kept below a certain threshold.
Since $x_6(t_f)$ and $x_7(t_f)$ depend on the random reaction rates and uncertain initial conditions, joint chance constraints are considered with required minimum constraint satisfaction probability $\beta=0.98$ and a confidence level of $1-\alpha = 0.99$.
It can be seen in Fig.~\ref{fig:CC:minSamples} that the desired satisfaction probability is guaranteed for $\beta_\text{cor} = 0.985$ and a sample size of $N_S = 5\,000$ samples. This implies
\begin{align}
\Prob{
\arraycolsep=1pt
\begin{array}{rcll}
x_6(t_{F}) \leq 0.6 \\
x_7(t_{F}) \leq 7.0
\end{array}
} \ge \beta_\text{cor} = 0.985. \label{eq:EXAMPLE:CHANCE_CONSTRAINT}
\end{align}

\subsubsection{Polynomial Chaos Expansion}

The system \eqref{eq:WOR} is not polynomial, which makes computation of the integrals in the Galerkin projections difficult (\cf Sec.~\ref{sec:PCE:Galerkin}). 
A polynomial-in-the-states representation was obtained by defining a new state $x_8(t) \coloneqq x_7(t)^{-1}$ and its derivative by $\dot{x}_8(t) = -u_1(t)/x_7(t)^2 = -x_8(t)^2u_1(t)$.
With this reformulation, \eqref{eq:WOR} becomes polynomial-in-the-states where each term is separable in the states, parameters and the inputs.
However, the employed transformations of the initial conditions $x_8(t_k) = x_7(t_k)^{-1}$ are not polynomials.
Hence, the collocation approach in Sec.~\ref{sec:PCE:collocation} was exploited to obtain a PC approximation (of order $P= 3$) of $x_8(t_k)$ by sampling $x_7(t_k)$. From the PC approximation of $x_8(t_k)$, the coefficients $\tbs{x_8}(t_k)$ were used as the initial conditions for the extended set of ordinary differential equations obtained by Galerkin projections \eqref{eq:sys:PCE}.

Note that similar transformations are always possible if the system has analytic nonlinearities \cite{Ohtsuka_2005_IEEETAC_Immersion} such as exponential terms, rational functions, etc.
Therefore, the Galerkin-based PC expansion as presented in Sec.~\ref{sec:PC} is applicable to a broad class of systems.

The PC expansion in conjunction with the Galerkin projection was applied to the polynomial-in-the-states system using a PC order $P=3$. The random variables were the reaction rates $k_1$, $k_2$, $k_3$ and the initial conditions $x_1(t_k), \ldots, x_7(t_k)$.
This lead to $\widetilde{P}=286$ terms in the expansion of each state (\cf Eq.~\eqref{eq:PCE:trunc}), resulting in 2\,002 differential equations altogether.
The first-order sensitivities of the states with respect to the inputs were used to compute the sensitivities and gradients, which resulted in a sensitivity system of 27\,456 differential equations.

\subsection{Closed-Loop Simulations}
\label{sec:Example:closed-loop}

The shrinking horizon SNMPC was implemented with batch end-time $t_f = 4\,000 \text{\,s}$ and sampling times $t_k \in \{0,\,250,\,500\,\,\dots,\,3\,500,\,3\,750\}$.
A piecewise-constant input parametrization was chosen with the initial values $u_{1}(t) = 0.002$ and $u_{2}(t) = 318$.

In the closed-loop simulations, the nonlinear model \eqref{eq:WOR} was used as the true plant, whose random parameters were drawn from the uncertainty distributions given above.
To update the the controller's state information at the beginning of each sampling time interval, normally distributed measurement/observer noise with a standard deviation of 1\,\% of the mean value was considered.  

\begin{figure*}[t!]
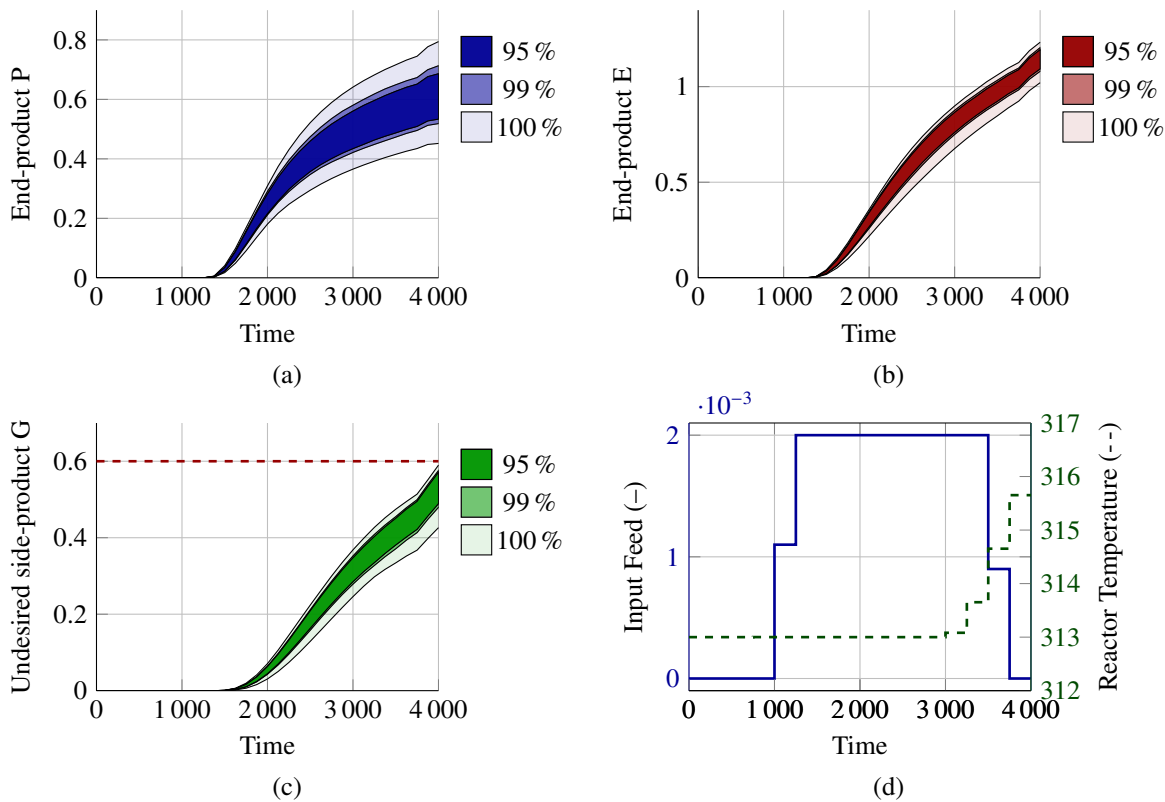

\centering
\begin{tabular}{cc}
\input{x4.tex} & \input{x5.tex} \\
(a) & (b)\\ 
\input{x6.tex} & \hspace{-0.5cm} 
%
%
\begin{tikzpicture}

\begin{axis}[%
width=4.5cm,
height=3.5491935483871cm,
scale only axis,
xmin=0,
xmax=4000,
xlabel={Time},
xmajorgrids,
separate axis lines,
every outer y axis line/.append style={blue},
every y tick label/.append style={font=\color{blue}},
ymin=-0.0001,
ymax=0.0021,
ytick={    0, 0.001, 0.002},
ylabel={Input Feed (--)},
ymajorgrids
]
\addplot[const plot,color=blue,solid,line width=1.0pt] plot table[row sep=crcr] {0	0\\
250	0\\
500	0\\
750	0\\
1000	0.0011\\
1250	0.002\\
1500	0.002\\
1750	0.002\\
2000	0.002\\
2250	0.002\\
2500	0.002\\
2750	0.002\\
3000	0.002\\
3250	0.002\\
3500	0.0009\\
3750	0\\
4000	0\\
};
\end{axis}

\begin{axis}[%
width=4.5cm,
height=3.5491935483871cm,
scale only axis,
xmin=0,
xmax=4000,
every outer y axis line/.append style={black!50!green},
every y tick label/.append style={font=\color{black!50!green}},
ymin=312,
ymax=317,
ytick={312, 313, 314, 315, 316, 317},
ylabel={Reactor Temperature (-\,-)},
axis x line*=bottom,
axis y line*=right
]
\addplot[const plot,color=black!50!green,dashed,line width=1.0pt] plot table[row sep=crcr] {0	313\\
250	313\\
500	313\\
750	313.001\\
1000	313\\
1250	313\\
1500	313\\
1750	313\\
2000	313\\
2250	313\\
2500	313\\
2750	313\\
3000	313.0814\\
3250	313.6521\\
3500	314.6521\\
3750	315.6521\\
4000	315.6521\\
};
\end{axis}
\end{tikzpicture}%
\\
(c) & (d)
\end{tabular}
\caption{Shrinking horizon {SNMPC}.
(a), (b) and (c) Concentrations of the end-products $P$ ($x_4$), $E$ ($x_5$), and side-product $G$ ($x_6$). Depicted are regions containing 90\,\%, 95\,\% and 100\,\% of 350 closed-loop simulations. As can be seen in (c), the uncertainties lead to constraint violations in no cases. (d) Optimal input profile for one representative sample.}
\label{fig:WOR:uncertain}
\end{figure*}

\subsection{Numerical Results}
\label{sec:Example:results}

\begin{figure*}[t!]
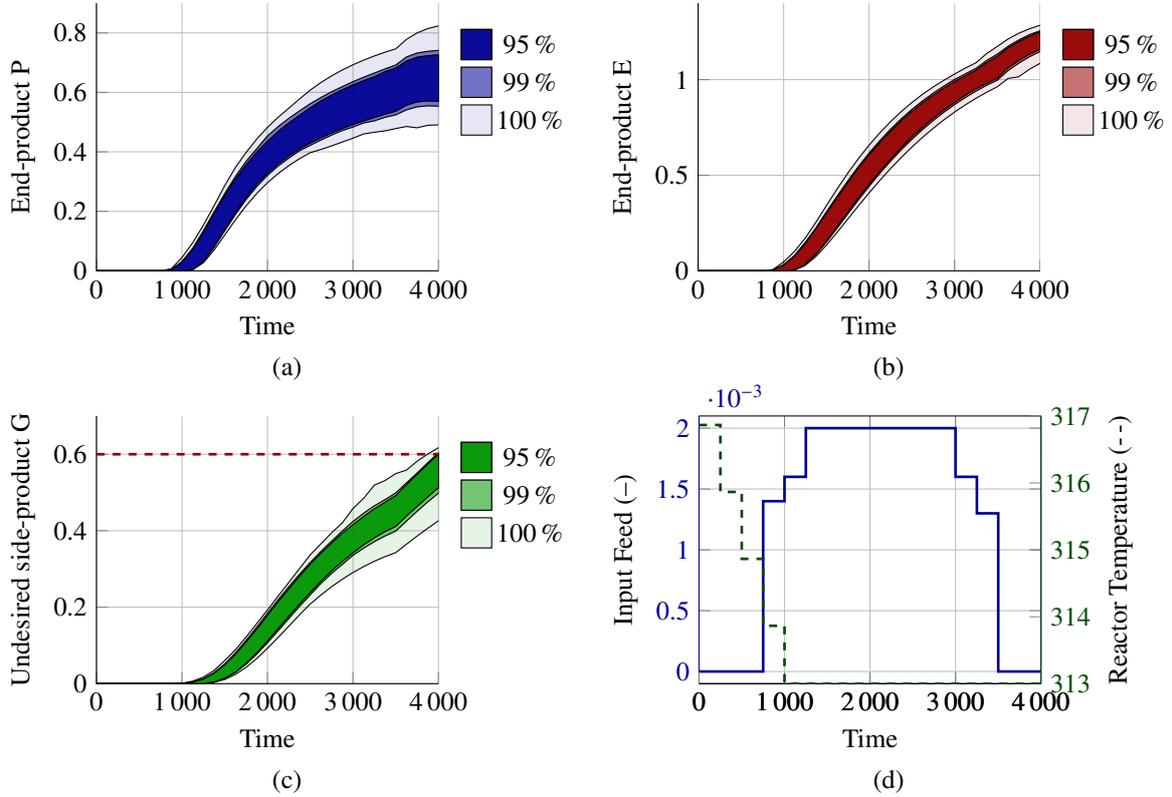

\centering
\begin{tabular}{cc}
\input{x4_nom.tex} & \input{x5_nom.tex} \\
(a) & (b)\\ 
\input{x6_nom.tex} & \hspace{-0.5cm} 
%
%
\begin{tikzpicture}

\begin{axis}[%
width=4.5cm,
height=3.5491935483871cm,
scale only axis,
xmin=0,
xmax=4000,
xlabel={Time},
xmajorgrids,
separate axis lines,
every outer y axis line/.append style={blue},
every y tick label/.append style={font=\color{blue}},
ymin=-0.0001,
ymax=0.0021,
ytick={     0, 0.0005,  0.001, 0.0015,  0.002},
ylabel={Input Feed (--)},
ymajorgrids
]
\addplot[const plot,color=blue,solid,line width=1.0pt] plot table[row sep=crcr] {0	0\\
250	0\\
500	0\\
750	0.0014\\
1000	0.0016\\
1250	0.002\\
1500	0.002\\
1750	0.002\\
2000	0.002\\
2250	0.002\\
2500	0.002\\
2750	0.002\\
3000	0.0016\\
3250	0.0013\\
3500	0\\
3750	0\\
4000	0\\
};
\end{axis}

\begin{axis}[%
width=4.5cm,
height=3.5491935483871cm,
scale only axis,
xmin=0,
xmax=4000,
every outer y axis line/.append style={black!50!green},
every y tick label/.append style={font=\color{black!50!green}},
ymin=313,
ymax=317,
ytick={313, 314, 315, 316, 317},
ylabel={Reactor Temperature (-\,-)},
axis x line*=bottom,
axis y line*=right
]
\addplot[const plot,color=black!50!green,dashed,line width=1.0pt] plot table[row sep=crcr] {0	316.8655\\
250	315.8653\\
500	314.8651\\
750	313.8648\\
1000	313\\
1250	313\\
1500	313\\
1750	313\\
2000	313\\
2250	313\\
2500	313\\
2750	313\\
3000	313\\
3250	313\\
3500	313\\
3750	313\\
4000	313\\
};
\end{axis}
\end{tikzpicture}%
\\
(c) & (d)
\end{tabular}
\caption{Shrinking horizon {NMPC} obtained by considering the mean values for each uncertain variable.
(a), (b) and (c) Concentrations of the end-products $P$ ($x_4$), $E$ ($x_5$), and side-product $G$ ($x_6$). Depicted are regions containing 90\,\%, 95\,\% and 100\,\% of 1\,000 closed-loop simulations. Note that in 225 cases (22.5\,\%) \texttt{fmincon} failed to find a feasible solution. Out of the remaining 775 cases, 37 (4.8\,\%) do not satisfy the nominal MPC's inequality constraint $x_6(t_f) \leq 0.6$.
(d) Optimal input profile for one representative sample. }
\label{fig:WOR:nominal}
\end{figure*}

The SNMPC was implemented in Matlab R2012a using \texttt{fmincon} with SQP method.
Time-critical code (including the ODEs, the Jacobians, the sensitivity equations) was written in C using Matlab CMEX-functions.
In particular, the integration of the differential equations was done using the \textit{SUNDIALS} integrator suite \cite{Hindmarsh_etAl_2005__SUNDIALS} using the nonlinear integrator \texttt{Functional} in conjunction with the \texttt{Adams} method.
The optimization was run on a Linux PC (Intel\texttrademark Core\texttrademark 2 Q6600, 2.4 GHz, 4 cores, 4 GB Ram).

\begin{table}[h]
\centering
\begin{tabular}{lr}
\hline
\multicolumn{2}{l}{\cellcolor[gray]{0.8}\textbf{Computing Times for Stochastic NMPC}} \\
\hline
\multicolumn{2}{l}{\textbf{Joint chance constraints} (5\,000 samples, 12 calls)}\\
\quad Probability of satisfaction and gradients & 2.3\,s \\
\hline
\textbf{Integration} (24 calls)  \\
\quad System dynamics & 65\,s\\
\quad Evaluation of the Jacobian & 162\,s \\
\quad Sensitivity differential equations & 288\,s \\
\hline
\multicolumn{2}{l}{\textbf{Total time} (including \texttt{fmincon} and }\\
\multicolumn{1}{l}{\quad further auxiliary functions)} & 404\,s\\
\hline
\end{tabular}
\caption{Averaged computing times in seconds for the chance constrained optimal control problem for the initial time-horizon $[0,4\,000]$.
The optimization times for the subsequent sampling time-points were much shorter and the total time ranged from $60$ to $150$ seconds.}
\label{tab:computationTime}
\end{table}

The results of the closed-loop SNMPC are shown in Fig.~\ref{fig:WOR:uncertain} and, for comparison, the results from the closed-loop NMPC are depicted in Fig.~\ref{fig:WOR:nominal}.
The results show that the presented SNMPC approach guarantees satisfaction of the joint chance constraint.
In contrast, the nominal NMPC caused violation of the constraints (\ie infeasible solutions) in about $25\%$ of the cases.

Table~\ref{tab:computationTime} summarizes the optimization time of the chance constrained optimal control problem at the first sampling time-point.
It can be seen that the evaluation of the Jacobians and sensitivity equations, which are required for the computation of the gradients, took the longest.
In case of finite-differencing approximation of the gradients, the optimization of the optimal control problem of the first horizon was prematurely stopped after several hours without finding a solution.
This demonstrates that the presented approximation of the gradients speeds up the optimization significantly. 

\begin{figure}[t!]
\centering
%
%
\begin{tikzpicture}

\begin{axis}[%
clip=false,
width=4.5cm,
height=3.5491935483871cm,
scale only axis,
xmin=2000,
xmax=4000,
xlabel={Time},
xmajorgrids=true,
xminorgrids=true,
ymode=log,
ymin=1e-8,
ymax=1e0,
yminorticks=true,
ylabel={Absolute Error},
ymajorgrids=true,
yminorgrids=true,
extra y ticks={1e-8,1e-6,1e-2,1e0}
]
\node[anchor=west] at (axis cs:4000,0.002127196428387) {\textcolor{blue}{mean}};
\node[anchor=west] at (axis cs:4000,5.9939992625926e-05) {\textcolor{red}{variance}};
\node[anchor=west] at (axis cs:4000,1.38914474725919e-06) {\textcolor{black}{kurtosis}};
\node[anchor=west] at (axis cs:4000,1.0315883172851e-05) {\textcolor{green}{skewness}};
\addplot [color=blue,solid,forget plot,line width=1.0pt]
  table[row sep=crcr]{
2000	0.000604172667927505\\
2125	0.000761432567893974\\
2250	0.000913367019077432\\
2375	0.00107273623607698\\
2500	0.00119247678521783\\
2625	0.00134436222084761\\
2750	0.00147151349198765\\
2875	0.00158974799170392\\
3000	0.00170758339336269\\
3125	0.00178006996214147\\
3250	0.00187545036218317\\
3375	0.00190039730508007\\
3500	0.0019686125108942\\
3625	0.00195610346990305\\
3750	0.00199615285168608\\
3875	0.00206587113967566\\
4000	0.00212719642838721\\
};
\addplot [color=red,dashed,forget plot,line width=1.0pt]
  table[row sep=crcr]{
2000	1.52224942502372e-05\\
2125	2.5769827810668e-05\\
2250	3.67568514495e-05\\
2375	4.65689271437483e-05\\
2500	5.41395596854073e-05\\
2625	5.89136552433562e-05\\
2750	6.07497343211138e-05\\
2875	6.04809192125554e-05\\
3000	5.79252816874719e-05\\
3125	5.45454516573788e-05\\
3250	4.98653058805198e-05\\
3375	4.54797016395904e-05\\
3500	4.05373320761411e-05\\
3625	3.66036835483556e-05\\
3750	3.25291066728501e-05\\
3875	3.10813973228206e-05\\
4000	2.9939992625926e-05\\
};
\addplot [color=green,dash pattern=on 1pt off 3pt on 3pt off 3pt,forget plot,line width=1.0pt]
  table[row sep=crcr]{
2000	4.39421662915867e-07\\
2125	1.03913766348908e-06\\
2250	1.9660384186706e-06\\
2375	3.16292527866287e-06\\
2500	4.52693359326125e-06\\
2625	5.92365186657154e-06\\
2750	7.26918934699435e-06\\
2875	8.4952059297287e-06\\
3000	9.57807307367698e-06\\
3125	1.05231414157542e-05\\
3250	1.13236942515061e-05\\
3375	1.20203180518251e-05\\
3500	1.25833053653996e-05\\
3625	1.30635702573628e-05\\
3750	1.34017291728947e-05\\
3875	1.4908518973255e-05\\
4000	1.6315883172851e-05\\
};
\addplot [color=black,dotted,forget plot,line width=1.0pt]
  table[row sep=crcr]{
2000	1.45767612928027e-08\\
2125	4.57651229210645e-08\\
2250	9.94987839474129e-08\\
2375	1.69100817011464e-07\\
2500	2.4182299962081e-07\\
2625	3.04320830467602e-07\\
2750	3.49464042554077e-07\\
2875	3.82272152451321e-07\\
3000	4.04305535930081e-07\\
3125	4.35238907976381e-07\\
3250	4.72534839816442e-07\\
3375	5.38358777642849e-07\\
3500	6.16473955582275e-07\\
3625	7.26736599700437e-07\\
3750	8.40293009745748e-07\\
3875	1.09173634827884e-06\\
4000	1.38914474725919e-06\\
};
\end{axis}
\end{tikzpicture}%
\caption{Accuracy of the PCE approximation illustrated through a comparison of 5\,000 Monte Carlo samples and samples generated using the PCE. Depicted is the absolute error of the mean, the variance, the skewness and kurtosis for the concentration of $G$ ($x_6(t))$.
Note that the curves are not plotted for $t<2\,000$ since the values are very small (i.\,e.\ approach $-\infty$ in the semi-logarithmic plot).
}
\label{fig:COMP:MC:PCE}
\end{figure}
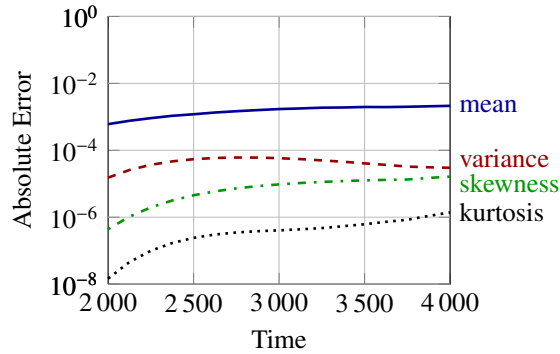

The simulation and evaluation of $5\,000$ samples using the PC approach was on average about 100 times faster (about 2\,s) than simulation of the same number of Monte-Carlo samples (about 200\,s).
This also emphasizes the advantage of the PC approach over a pure sampling or scenario-based approaches in the investigated context. The accuracy of the PC approximation is compared with Monte-Carlo samples in Fig.~\ref{fig:COMP:MC:PCE} demonstrating adequate  approximation quality.

\section{Conclusions and Outlook}
\label{sec:Conclusions}

This work demonstrates the feasibility of a SNMPC approach for real-time control of a chemical process with uncertain parameters and initial conditions.
Polynomial chaos is used for uncertainty propagation and Monte Carlo sampling, which is significantly more efficient than Monte Carlo simulations based on the nonlinear system model.
The approximation of the probability densities reduces to solving an expanded set of differential equations (to get the values of the PC coefficients) and matrix multiplications to evaluate the chance constraints.
The proposed sample approximation of the joint chance constraints is very flexible and the gradient approximation improves efficiency of the overall optimization significantly.
The confidence analysis provides a systematic way to determine the sample size and the chance constraint tightening required to guarantee a desired feasibility probability and constraint satisfaction probability.
This allows the systematic trade-off between large sample sizes (\ie accurate results, time-consuming) and small sample sizes (\ie less accurate results, less time-consuming).
This could also be useful when one wants to adapt the \mbox{SNMPC} approach to hardware with less computational power, or where timing-constraints have to be satisfied.
Note that the sample approximations and confidence analysis is quite general and independent of the used method to generate samples, and independent of the properties (\eg convexity) of the optimization problem.

It is important to remark that the probabilities of satisfaction and feasibility hold only for the PC approximation, but not necessarily for the original uncertain nonlinear system.
However, as illustrated in the example in this work and mentioned in other work (see introduction), PC allows for accurate predictions of the propagation of stochastic uncertainties through (nonlinear) system dynamics.
However, the accuracy is clearly system dependent and to the best of our knowledge no systematic method exist to determine a-priori or depending on the system dynamics the PC order required to achieve a certain accuracy of the estimates.
However, there has been recent progress in the error analysis of PC expansions \cite{Shi_Zhang_2012_ErrorAnalys_PCE} and a line of future research could be to include the error analysis in the synthesis of a robust controller.

From Eq.~\eqref{eq:PCE} it becomes clear that the number of monomials required in the PCE grows rapidly with the chosen PC order and number of uncertain variables.
Using sparse PCE, \ie different polynomial orders for different variables \cite{Luo_2006_PhD_PCE_PDE}, is one solution to reduce this complexity.
It has been noted \cite{Gerritsma_etAl_2010_TD_PCE} that PCE provides accurate predictions for small time only and that the predictions may become inaccurate for increasing time.
We believe, however, that the error is negligible in a controlled system with constant update by measurements.
In any case, the prediction could be improved by adapting the orthogonal basis as proposed in \cite{Gerritsma_etAl_2010_TD_PCE}.
Another interesting future extension could be the consideration of time-dependent disturbances (\eg \cite{Huschto_Sager_2013_ECC__OC_PCE}).

\section*{Acknowledgements}

The authors thank Rolf Findeisen for support, and Felix Petzke for helping with the implementation of the example.
We also thank Bhushan Gopaluni, Timm Faulwasser and Philipp Rumschinski for their helpful comments improving the manuscript.

\appendix
\section{Derivation and Proof of Theorem~1}

To determine the least conservative constraint tightening, a statistical analysis for the sample approximation \eqref{eq:prob:approx} (with a given number of samples $N_S$ and confidence level $\alpha$) is used.

Let $p \coloneqq \Prob{\bs{g}\left(\tbs{X}, \bs{\xi}\right) \le \bs{0} }$ (cf.\ \eqref{eq:prob}), and the sample approximation $\widehat{p} \coloneqq \frac{1}{N_S}
\sum_{j=1}^{N_S}I_{\mathcal{G}}\left(\bs{\xi}^{[j]}\right)$ (cf.\ \eqref{eq:prob:approx}).
The estimation of $p$ is a well-known problem in statistics and corresponds to estimating the success probability of a sequence of Bernuoulli trials \cite{Hollander_etAl_2013_Nonparametric_Statistical_Methods,Rossi_CC_ConfidenceInterval,Stan_etAl_2014_CC_ConfidenceInterval}.
This is due to the fact that testing whether a sample $\bs{\xi}^{[j]}$ lies in  $\mathcal{G}$ (i.\,e., evaluation of $I_\mathcal{G}(\bs{\xi}^{[j]})$) is a Bernoulli trial with a ``yes'' or ``no"  outcome (i.\,e., ``satisfied'' or ``violated'').
Therefore, $\widehat{p}$ is a realization of the random variable $\widehat{P}$ that is distributed according to a binomial distribution, $\widehat{P} \sim \frac{1}{N_S} \text{Bin}(N_S, p)$.
The confidence interval $[p_\text{lb}(\alpha, N_S, \widehat{P}),p_\text{ub}(\alpha, N_S, \widehat{P})]$ consists of a range of values that, with a certain confidence level of $1-\alpha$, contains the true value $p$.
Furthermore, the probability of $p \ge \beta$ given the observation that $p_\text{lb}(\alpha, N_S, \widehat{P}) \ge \beta$ is not less than $1-\alpha$:
\begin{align}
\label{eq:qlb}
\Prob{p\ge\beta\ |\ p_\text{lb}(\alpha, N_S, \widehat{P}) \ge \beta} \ge 1-\alpha.
\end{align}
Since $p_\text{lb}(\alpha, N_S, \widehat{P})$ is monotonically increasing in $\widehat{P}$, it holds that
\begin{align}
\label{eq:qlb:iff}
\exists\ \beta_\text{cor} \in [\beta,1]\text{\ s.\,t.\ } \widehat{P} \ge \beta_\text{cor}\text{\ iff\ }p_\text{lb}(\alpha, N_S, \widehat{P}) \ge \beta.
\end{align}
It follows from \eqref{eq:qlb} and \eqref{eq:qlb:iff} that
\begin{align*}
\Prob{p\ge\beta\ |\ \widehat{P} \ge \beta_\text{cor}} \ge 1-\alpha.
\end{align*}
Thus, the lower confidence bound is used to determine $\beta_\text{cor}$, since we demand that $p_\text{lb}(\alpha, N_S, \widehat{P}) \ge \beta$, meaning that the confidence interval shall not cover the infeasible interval $p<\beta$.
The lower confidence bound can be determined approximatively from the quantiles of a normal distribution or exactly from the percentiles of the inverse cumulative Beta distribution $\text{betainv}$ \cite{Hollander_etAl_2013_Nonparametric_Statistical_Methods,Rossi_CC_ConfidenceInterval,Stan_etAl_2014_CC_ConfidenceInterval}  as
\begin{equation}
\label{eq:alpha:Beta}
p_\text{lb}(\alpha, N_S, \floor{\beta_\text{cor} N_S}) =  1- \text{betainv}\left(1 - \frac{\alpha}{2},N_S +1 - \floor{\beta_\text{cor} N_S}, \floor{\beta_\text{cor} N_S}\right)\nonumber.
\end{equation}
Here $\floor{\beta_\text{cor} N_S}\in \mathbb{N}$ is the number of ``satisfied'' Bernoulli trials,
and $\beta_\text{cor}$ represents a corrected (i.\,e.\ tightened) satisfaction probability.
Since the quantiles are exact and due to the equivalence above, it follows that $\beta_\text{cor}$ is the least conservative constraint tightening.
The statement \eqref{eq:N_S} in the theorem follows directly from \eqref{eq:alpha:Beta} and from the requirement $p_\text{lb}(\alpha, N_S, \floor{\beta_\text{cor} N_S}) \ge \beta$. \hfill$\square$

\section{Proof of Theorem~2}

Sampling only affects the cost function \eqref{eq:prob1:cost} and the chance constraints \eqref{eq:prob1:cc}, which means that $\bs{\pi}'$ satisfies all other constraints in \eqref{eq:prob1}.
Due to Theorem~1, $\bs{\pi}'$ satisfies the chance constraint \eqref{eq:prob1:cc} with a confidence level not less than $1-\alpha$.
From this it can be concluded that $\bs{\pi}'$ is a feasible point of \eqref{eq:prob1} with a probability not less than $1-\alpha$.
$\hfill\square$

Note that $\bs{\pi}'$ is not necessarily an optimal solution to \eqref{eq:prob1} because the sample approximation of the cost function may differ from the exact value of the cost function.

\section{Derivative and Proof of Proposition~1}

For simpler presentation, we provide the proof for a scalar function $g(\cdot)$ (\ie for individual chance constraints) first.
The extension to vector-valued functions $\bs{g}(\cdot)$ (\ie joint chance constraints) is straightforward and is done last.

The proof follows from the construction as described in the following.

\subsection*{Reformulations}

Since the first-order sensitivities are assumed to be given, what remains to be computed are the partial derivatives $\frac{\partial \Prob{\cdot}}{\partial \tbs{X}}$
For this purpose, we derive another approximation of the satisfaction probability \eqref{eq:prob}, which then allows to determine $\frac{\partial \Prob{\cdot}}
{\partial \tbs{X}}$.
The main idea is to approximate the $n_\xi$-dimensional integration in \eqref{eq:prob} by integration over only one random variable, say $\xi_1$, while keeping the remaining $n_\xi-1$ variables in $\obs{\xi}\coloneqq [\xi_2, \ldots, \xi_{n_\xi} ]$ fixed to values as determined by sampling.
We assumed, without loss of generality, that we integrate with respect to $\xi_1$. The analysis can be performed analogously for any other random variable $\xi_i$, $i=2,\ldots,n_\xi$.

Assume (for now) that for each sample $\obs{\xi}^{[j]}$, a vector defining the integration limits $\bs{r}^{[j]} \coloneqq [r_0^{[j]}, r_1^{[j]}, \ldots, r_{n_r^{[j]}+1}^{[j]}]$ (see Fig.~\ref{fig:CC:approach}b) is given.
\footnote{The integration limits formally depend on $\tbs{X}$ and on the sample $\obs{\xi}^{[j]}$, and that the number $n_{r^{[j]}}$ may be different for each sample.
See Fig.~\ref{fig:CC:approach}b for further explanations.
However, we omit this dependency for shorter notation.}
With these integration limits given, Eq.~\eqref{eq:prob} can be approximated\footnote{Eq.~\eqref{eq:prob:h} provides a better approximation of \eqref{eq:prob} than \eqref{eq:prob:approx} since the entire probability distribution of $\xi_1$ is taken into account; however, this approximation requires the evaluation of the cumulative probability density function of $\mu_1$, which is for many distributions not available in an analytic form.}
  as 
\begin{equation}
\label{eq:prob:h}
\Prob{g\left(\tbs{X}, \xi_1,\obs{\xi}\right) \le \bs{0} }
\approx\\
\frac{1}{N_S}\sum_{j=1}^{N_S}
\left[
\sum_{i=0}^{n_{r^{[j]}}}
I_{\mathcal{G}}\left(\hat r^{[j]}_i,\obs{\xi}^{[j]}\right)
\lint_{r^{[j]}_{i}}^{r^{[j]}_{i+1}}
\mu_{1}(d\xi_1)
\right],
\end{equation}
where the indicator function is evaluated at the point $(\hat{r}_i^{[j]},\obs{\xi}^{[j]})$, with $\hat{r}^{[j]}_i$ being an arbitrary point from the interior of the integration interval $(r^{[j]}_{i},r^{[j]}_{i+1} )$, as e.\,g.\ the mid-point $\hat{r}_i^{[j]} = \frac{r_i^{[j]}+r_{i+1}^{[j]}}{2}$ in \eqref{eq:gradProb:h}.


Based on the vector of integration bounds $\bs{r}^{[j]}$, the gradient can be approximated by derivative of \eqref{eq:prob:h} with respect to $\tbs{X}$, which gives \eqref{eq:gradProb:h}.
Note that we used the facts that the order of integration and differentiation can be changed.
Since the integration limits $r_{i}^{[j]}$ in Eq.~\eqref{eq:prob:h} depend on the (known and constant) $\tbs{X}$, the partial derivatives $\pfrac{r_{i}^{[j]}}{\tbs{X}}$ have to be taken into account, and they are given by \eqref{eq:gradProb:h:implicitFunc}.
The first row in \eqref{eq:gradProb:h:implicitFunc} follows from the fact that $r_{0}^{[j]}$ and $r_{n_{r^{[j]}}+1}^{[j]}$ are constants and defined by the minimum and maximum of the support of $\xi_1$; the second row in \eqref{eq:gradProb:h:implicitFunc} follows from the implicit function theorem.

From Eq.~\eqref{eq:gradProb:h:implicitFunc} it becomes clear that one has to avoid division by zero, which is the case if $g\left(\tbs{X}, \xi_1, \obs{\xi}\right) = 0$ has multiple roots with respect to $\xi_1$.
These cases have to be treated and are first discussed for individual chance constraints. The extension to joint chance constraints follows immediately from that.

\subsection*{Integration Limits for Individual Chance Constraints}

As illustrated in Fig.~\ref{fig:CC:approach}, the integration limits are defined by the solutions of $g\left(\tbs{X}, \xi_1,\obs{\xi}^{[j]}\right) = 0$ or, in other words, the values of $\xi_1$ where the indicator function $I_\mathcal{G}\left(r_i^{[j]},\obs{\xi}^{[j]}\right)$ changes its value when traversing along the direction of $\xi_1$ for fixed values of $\obs{\xi}^{[j]}$.
In the first step, we therefore solve the equations $g_i\left(\tbs{X}, \xi_1,\obs{\xi}^{[j]}\right) = 0$ for $\xi_1$ and determine all its (real-valued) roots on the support of the random variable.
Due to the PC expansion involving polynomials in $\xi_1$ of order usually greater than one, the equation is possible to have no unique (real-valued) solution.
From Eq.~\eqref{eq:gradProb:h:implicitFunc} it is clear that roots with multiplicity greater than 1 lead to division by zero and therefore cannot be considered, and the sample has to be discarded.

In the next step, all (real-valued) roots are sorted such that $r_{1}^{[j]} \le r_{2}^{[j]} \le \cdots \le r_{n_{r^{[j]}}}^{[j]}$, where $n_{r^{[j]}}$ denotes the number of real-valued roots.
As the last step, we introduce $r_{0}^{[j]}$ (resp.\ $r_{n_{r^{[j]}}+1}^{[j]}$) corresponding to the lowest (resp.\ largest) value that $\xi_1$ can take on its support (which may also be $\pm \infty$).
With that, one obtains the vector $\bs{r}^{[j]} \coloneqq [r_0^{[j]}, r_1^{[j]}, \ldots, r_{n_{r^{[j]}}+1}^{[j]}]$.

\subsection*{Integration Limits for Joint Chance Constraints}

One proceeds similar as for individual chance constraints and determines the real-valued roots on the support of $\xi_1$ for each equation $g_i(\cdot)$, $i=1,\ldots,N_g$.
As for individual chance constraints, samples having (real-valued) roots with multiplicity greater than one are discarded.
In addition, one has to take care of the fact that roots of the different equations may coincide.
To avoid such degenerate cases and the ambiguity in applying equation \eqref{eq:gradProb:h:implicitFunc}, such samples are also discarded.
After that, all roots as well as $r_{0}^{[j]}$ and $r_{n_{r^{[j]}+1}}^{[j]}$ (see individual chance constraints) are collected in vector in ascending order.
When evaluating \eqref{eq:gradProb:h:implicitFunc}, one has to use the function $g_i(\cdot)$ for which $r_{i}^{[j]}$ is a root.

With that, all elements and equations are established, which completes the proof. \hfill$\square$

\biboptions{sort&compress}
\bibliographystyle{elsarticle-num}
\bibliography{PCE_SNMPC_CC_SA}

\begin{thebibliography}{10}
\expandafter\ifx\csname url\endcsname\relax
  \def\url#1{\texttt{#1}}\fi
\expandafter\ifx\csname urlprefix\endcsname\relax\def\urlprefix{URL }\fi
\expandafter\ifx\csname href\endcsname\relax
  \def\href#1#2{#2} \def\path#1{#1}\fi

\bibitem{Qin_etAl_1997_Overview_industrial_MPC}
S.~J. Qin, T.~A. Badgwell, An overview of industrial model predictive control
  technology, in: AIChE Symposium Series, Vol.~93, American Institute of
  Chemical Engineers, 1997, pp. 232--256.

\bibitem{Rawlings_2000_Tutorial_Overview_MPC}
J.~B. Rawlings, Tutorial overview of model predictive control, IEEE Control
  Systems Magazine 20~(3) (2000) 38--52.

\bibitem{mor99}
M.~Morari, J.~H. Lee, Model predictive control: past, present and future,
  Computers \& Chemical Engineering 23~(4-5) (1999) 667--682.

\bibitem{Bemporad_Morari_1999_RMPC_Survey}
A.~Bemporad, M.~Morari, Robust model predictive control: {A} survey, in:
  A.~Garulli, A.~Tesi (Eds.), Robustness in Identification and Control,
  Springer, Berlin, London, New York, 1999, pp. 207--226.

\bibitem{Garatti_Campi_2013_CSM__Modulating_Robustness_ControlDesign}
S.~Garatti, M.~C. Campi, Modulating robustness in control design: {P}rinciples
  and algorithms, IEEE Control Systems Magazine 33~(2) (2013) 36--51.

\bibitem{Vidyasagar_2001_Automatica__RA_Rob_ContrSyn}
M.~Vidyasagar, Randomized algorithms for robust controller synthesis using
  statistical learning theory, Automatica 37 (2001) 1515--1528.

\bibitem{sch99}
A.~Schwarm, M.~Nikolaou, Chance-constrained model predictive control, AIChE
  Journal 45 (1999) 1743--1752.

\bibitem{Li_etAl_2000_RMPC_CC}
P.~Li, M.~Wendt, G.~Wozny, Robust model predictive control under chance
  constraints, Computers \& Chemical Engineering 24~(2) (2000) 829--834.

\bibitem{hes03}
D.~H.~V. Hessem, O.~H. Bosgra, A full solution to the constrained stochastic
  closed-loop {MPC} problem via state and innovations feedback and its receding
  horizon implementation, in: Proc. 42nd IEEE Conference on Decision and
  Control (CDC), Maui, 2003, pp. 929--934.

\bibitem{Li_etAl_2008_CC_PC_UncSys}
P.~Li, H.~Arellano-Garcia, G.~Wozny, Chance constrained programming approach to
  process optimization under uncertainty, Computers \& Chemical Engineering
  32~(1) (2008) 25--45.

\bibitem{Cannon_etAl_2011_Stochastic_TubeMPC_CC}
M.~Cannon, B.~Kouvaritakis, S.~V. Rakovic, Q.~Cheng, Stochastic tubes in model
  predictive control with probabilistic constraints, IEEE Transactions on
  Automatic Control 56~(1) (2011) 194--200.

\bibitem{pri09}
J.~Primbs, C.~Sung, Stochastic receding horizon control of constrained linear
  systems with state and control multiplicative noise, IEEE Transactions on
  Automatic Control 54 (2009) 221--230.

\bibitem{ber09}
D.~Bernardini, A.~Bemporad, Scenario-based model predictive control of
  stochastic constrained linear systems, in: Proc. 48th IEEE Conference on
  Decision and Control (CDC), Shanghai, 2009, pp. 6333--6338.

\bibitem{Cannon_etAl_2009_Automatica__MPC_StochMultUnc}
M.~Cannon, B.~Kouvaritakis, X.~Wu, Model predictive control for systems with
  stochastic multiplicative uncertainty and probabilistic constraints,
  Automatica 45~(1) (2009) 167 -- 172.

\bibitem{Zhang_etAl_2013_CDC__LinSMPC_RA_RO}
X.~Zhang, K.~Margellos, P.~Goulart, J.~Lygeros, Stochastic model predictive
  control using a combination of randomized and robust optimization, in: Proc.
  IEEE 52nd Annual Conference on Decision and Control (CDC), 2013, pp.
  7740--7745.

\bibitem{Blackmore_etAl_2010_IEEETRob_Particle_SMPC_CC}
L.~Blackmore, M.~Ono, A.~Bektassov, B.~C. Williams, A probabilistic
  particle-control approximation of chance-constrained stochastic predictive
  control, Robotics, IEEE Transactions on 26~(3) (2010) 502--517.

\bibitem{Oldewurtel_etAl_2008_CDC_Approx_CC_affineDisturbFB}
F.~Oldewurtel, C.~N. Jones, M.~Morari, A tractable approximation of chance
  constrained stochastic {MPC} based on affine disturbance feedback, in: Proc.
  47th IEEE Conference on Decision and Control (CDC), 2008, pp. 4731--4736.

\bibitem{Kouvaritakis_etAl_2013_IJSS_DisturbComp_SMPC}
B.~Kouvaritakis, M.~Cannon, D.~Mu{\~n}oz-Carpintero, Efficient prediction
  strategies for disturbance compensation in stochastic {MPC}, International
  Journal of Systems Science 44~(7) (2013) 1344--1353.

\bibitem{Korda_etAl_2014_IEEETAC__LinSMPC_AverageConstrViol}
M.~Korda, R.~Gondhalekar, F.~Oldewurtel, C.~N. Jones, Stochastic {MPC}
  framework for controlling the average constraint violation, IEEE Transactions
  on Automatic Control 59~(7) (2014) 1706--1721.

\bibitem{Hashimoto_2013_CDC__LinSMPC_Chebychev}
T.~Hashimoto, Probabilistic constrained model predictive control for linear
  discrete-time systems with additive stochastic disturbances, in: Proc. IEEE
  52nd Annual Conference on Decision and Control (CDC), 2013, pp. 6434--6439.

\bibitem{Farina_etAl_2013_CDC__LinSMPC_Cantelli}
M.~Farina, L.~Giulioni, L.~Magni, R.~Scattolini, A probabilistic approach to
  model predictive control, in: Proc. IEEE 52nd Annual Conference on Decision
  and Control (CDC), 2013, pp. 7734--7739.

\bibitem{Calafiore_Ghaoui_20016_JOTA__DistrRobust_CC_LP}
G.~C. Calafiore, L.~E. Ghaoui, On distributionally robust chance-constrained
  linear programs, Journal of Optimization Theory and Application 130~(1)
  (2006) 1--22.

\bibitem{Geletu_etAl_2013_IJSS__Review_CCOpt}
A.~Geletu, M.~Kl{\"o}ppel, H.~Zhang, P.~Li, Advances and applications of
  chance-constrained approaches to systems optimisation under uncertainty,
  International Journal of Systems Science 44~(7) (2013) 1209--1232.

\bibitem{Royset_Polak_2004_SAA_CC}
J.~O. Royset, E.~Polak, Reliability-based optimal design using sample average
  approximations, Probabilistic Engineering Mechanics 19~(4) (2004) 331--343.

\bibitem{Shapiro_2008_MP__StochProgr_UncOpt}
A.~Shapiro, Stochastic programming approach to optimization under uncertainty,
  Mathematical Programming 112 (2008) 183--220.

\bibitem{Calafiore_Campi_2006_IEEETAC__Scenario_Approach_RobustControl}
G.~C. Calafiore, M.~C. Campi, The scenario approach to robust control design,
  IEEE Transactions on Automatic Control 51~(5) (2006) 742--753.

\bibitem{Campi_Garatti_2008_SIAMJO__ExactFeasibility_RA_ConvexPrograms}
M.~C. Campi, S.~Garatti, The exact feasibility of randomized solutions of
  uncertain convex programs, SIAM Journal on Optimization 19 (2008) 1211--1230.

\bibitem{Calafiore_2010_SIAMJO__Random_ConvexPrograms}
G.~Calafiore, Random convex programs, SIAM Journal on Optimization 20~(6)
  (2010) 3427--3464.

\bibitem{Schildbach_etAl_2012_ACC__RA_LinMPC}
G.~Schildbach, G.~C. Calafiore, L.~Fagiano, M.~Morari, Randomized model
  predictive control for stochastic linear systems, in: Proc. American Control
  Conference (ACC), Montreal, 2012, pp. 417--422.

\bibitem{Calafiore_Fabiano_2013_IEEETAC__RMPC_ScenarioOpt}
G.~C. Calafiore, L.~Fagiano, Robust model predictive control via scenario
  optimization, IEEE Transactions on Automatic Control 58~(1) (2013) 219--224.

\bibitem{Alamo_etAl_2010_Sample_Complexity_Prob_Analysis_Design}
T.~Alamo, R.~Tempo, A.~Luque, On the sample complexity of probabilistic
  analysis and design methods, in: Perspectives in Mathematical System Theory,
  Control, and Signal Processing, Springer-Verlag London, London, 2010, pp.
  39--50.

\bibitem{Alamo_etAl_2014_RA_UncSys_Sample_Complexity}
T.~Alamo, R.~Tempo, A.~Luque, D.~Ramirez, Randomized methods for design of
  uncertain systems: Sample complexity and sequential algorithms, {arXiv}:
  1110.1892v2.

\bibitem{Geletu_etAl_2014_EOpt__AnalyticApprox_NonConvex_CC}
A.~Geletu, M.~Kl{\"o}ppel, A.~Hoffmann, P.~Li, A tractable approximation of
  non-convex chance constrained optimization with non-{G}aussian uncertainties,
  Engineering Optimization~(ahead-of-print) (2014) 1--26.

\bibitem{Feng_etAl_2011_CDC_KinshipFuncApprox_CC}
C.~Feng, F.~Dabbene, C.~M. Lagoa, A kinship function approach to robust and
  probabilistic optimization under polynomial uncertainty, IEEE Transactions on
  Automatic Control 56~(7) (2011) 1509--1523.

\bibitem{Mesbah_etAl_2014_ACC__SNMPC_CC}
A.~Mesbah, S.~Streif, R.~Findeisen, R.~D. Braatz, Stochastic nonlinear model
  predictive control with probabilistic constraints, in: Proc. American Control
  Conference (ACC), Portland, Oregon, 2014, pp. 2413--2419.

\bibitem{Imsland_etAl_2010_MB_Contr_Comparison_FiniteDiff_SensODE}
L.~Imsland, P.~Kittilsen, T.~S. Schei, Model-based optimizing control and
  estimation using modelica model, Modeling, Identification and Control 31~(3)
  (2010) 107--121.

\bibitem{Garnier_etAl_2009_Approximative_Gradients_CC}
J.~Garnier, A.~Omrane, Y.~Rouchdy, Asymptotic formulas for the derivatives of
  probability functions and their monte carlo estimations, European Journal of
  Operational Research 198~(3) (2009) 848--858.

\bibitem{Paulson_etAl_2014_Fast_SMPC}
J.~A. Paulson, A.~Mesbah, S.~Streif, R.~Findeisen, R.~D. Braatz, Fast
  stochastic model predictive control of high-dimensional systems, Vol. In
  Press, 2014.

\bibitem{fag12}
L.~Fagiano, M.~Khammash, Nonlinear stochastic model predictive control via
  regularized polynomial chaos expansions, in: Proc. 51st IEEE Conference on
  Decision and Control (CDC), Maui, 2012, pp. 142--147.

\bibitem{Huschto_Sager_2013_ECC__OC_PCE}
T.~Huschto, S.~Sager, Stochastic optimal control in the perspective of the
  wiener chaos, in: Proc. European Control Conference (ECC), Zurich, 2013, pp.
  3059--3064.

\bibitem{Kwang-Ki_Braatz_2013_IJC_MPC_PCE}
K.-K.~K. Kim, R.~D. Braatz, Generalised polynomial chaos expansion approaches
  to approximate stochastic model predictive control, International Journal of
  Control 86~(8) (2013) 1324--1337.

\bibitem{wie38}
N.~Wiener, The homogeneous chaos, American Journal of Mathematics 60 (1938)
  897--936.

\bibitem{Ghanem_Spanos_1991_SFE}
R.~Ghanem, P.~Spanos, Stochastic Finite Elements - A Spectral Approach,
  Springer-Verlag, New York, 1991.

\bibitem{xiu02}
D.~Xiu, G.~E. Karniadakis, The wiener-askey polynomial chaos for stochastic
  differential equations, SIAM Journal of Scientific Computation 24 (2002)
  619--644.

\bibitem{Kim_etal_2013_CSM__PCE_review}
K.-K.~K. Kim, D.~E. Shen, Z.~K. Nagy, R.~D. Braatz, Wiener's polynomial chaos
  for the analysis and control of nonlinear dynamical systems with
  probabilistic uncertainties, IEEE Control Systems Magazine 33~(5) (2013)
  58--67.

\bibitem{Ohtsuka_2005_IEEETAC_Immersion}
T.~Ohtsuka, Model structure simplification of nonlinear systems via immersion,
  IEEE Transactions on Automatic Control 50~(5) (2005) 607--618.

\bibitem{Fisher_Bhattacharya_2009_LQR_PCE}
J.~Fisher, R.~Bhattacharya, Linear quadratic regulation of systems with
  stochastic parameter uncertainties, Automatica 45~(12) (2009) 2831--2841.

\bibitem{Calafiore_Campi_2005_MathProgr_UncConvProg_RA_ConfidenceLevels}
G.~Calafiore, M.~C. Campi, Uncertain convex programs: randomized solutions and
  confidence levels, Mathematical Programming 102~(1) (2005) 25--46.

\bibitem{Cameron_Martin_1947_Orthogonal_NonLin_Func}
R.~H. Cameron, W.~T. Martin, The orthogonal development of non-linear
  functionals in series of {F}ourier-{H}ermite functionals, Annals of
  Mathematics 48 (1947) 385--392.

\bibitem{Oladyshkin_Nowak_2012_Arbitrary_PCE}
S.~Oladyshkin, W.~Nowak, Data-driven uncertainty quantification using the
  arbitrary polynomial chaos expansion, Reliability Engineering \& System
  Safety 106 (2012) 179--190.

\bibitem{Gerritsma_etAl_2010_TD_PCE}
M.~Gerritsma, J.-B. Van~der Steen, P.~Vos, G.~Karniadakis, Time-dependent
  generalized polynomial chaos, Journal of Computational Physics 229~(22)
  (2010) 8333--8363.

\bibitem{Nagy_Braatz_2007_JPC_UncertAnalysis_PCE}
Z.~Nagy, R.~Braatz, Distributional uncertainty analysis using power series and
  polynomial chaos expansions, Journal of Process Control 17~(3) (2007)
  229--240.

\bibitem{Mesbah_etAl_2014_IFACWC__aFDI_ProbUnc}
A.~Mesbah, S.~Streif, R.~Findeisen, R.~D. Braatz, Active fault diagnosis for
  nonlinear systems with probabilistic uncertainties, in: Proc. 19th IFAC World
  Congress, Cape Town, South Africa, 2014, pp. 7079--7084.

\bibitem{Tatang_etAl_1997_Efficient_ParamUnc_Collocation}
M.~A. Tatang, W.~Pan, R.~G. Prinn, G.~J. McRae, An efficient method for
  parametric uncertainty analysis of numerical geophysical models, Journal of
  Geophysical Research: Atmospheres (1984--2012) 102~(D18) (1997) 21925--21932.

\bibitem{Gautschi_etAl_2004_Book__Orthogonal_Polynomials}
W.~Gautschi, R.~S. Friedman, J.~Burns, R.~Darjee, A.~Mcintosh, Orthogonal
  Polynomials: Computation and Approximation, Numerical Mathematics and
  Scientific Computation Series, Oxford University Press, Oxford, U.K., 2004.

\bibitem{Streif_etAl_2014_IFACWC__PCE_ExpDesign_MD_PC}
S.~Streif, F.~Petzke, A.~Mesbah, R.~Findeisen, R.~D. Braatz, Optimal
  experimental design for probabilistic model discrimination using polynomial
  chaos, in: Proc. 19th IFAC World Congress, Cape Town, South Africa, 2014, pp.
  4103--4109.

\bibitem{Campi_Garatti_2010_JOTA__SA_CC}
M.~Campi, S.~Garatti, A sampling-and-discarding approach to chance-constrained
  optimization: Feasibility and optimality, Journal of Optimization Theory and
  Applications 148~(2) (2011) 257--280.

\bibitem{Rossi_CC_ConfidenceInterval}
R.~Rossi, Solving stochastic constraint programs via sampling, {arXiv}:
  1110.1892v2.

\bibitem{Hannemann_Marquardt_2010_Adjoints_Hessians_Lagragian_DynOpt_WilliamsOttoReactor}
R.~Hannemann, W.~Marquardt, Continuous and discrete composite adjoints for the
  {H}essian of the {L}agrangian in shooting algorithms for dynamic
  optimization, SIAM Journal on Scientific Computing 31~(6) (2010) 4675--4695.

\bibitem{Williams_Otto_1960_Reactor}
T.~J. Williams, R.~E. Otto, A generalized chemical processing model for the
  investigation of computer control, Transactions of the American Institute of
  Electrical Engineers, Part I: Communication and Electronics 79~(5) (1960)
  458--473.

\bibitem{Hindmarsh_etAl_2005__SUNDIALS}
A.~C. Hindmarsh, P.~N. Brown, K.~E. Grant, S.~L. Lee, R.~Serban, D.~E.
  Shumaker, C.~S. Woodward,
  \href{http://dx.doi.org/10.1145/1089014.1089020}{{SUNDIALS}: Suite of
  nonlinear and differential/algebraic equation solvers}, ACM Trans. Math.
  Softw. 31~(3) (2005) 363--396.
\newline\urlprefix\url{http://dx.doi.org/10.1145/1089014.1089020}

\bibitem{Shi_Zhang_2012_ErrorAnalys_PCE}
W.~Shi, C.~Zhang, Error analysis of generalized polynomial chaos for nonlinear
  random ordinary differential equations, Applied Numerical Mathematics 62~(12)
  (2012) 1954--1964.

\bibitem{Luo_2006_PhD_PCE_PDE}
W.~Luo, Wiener chaos expansion and numerical solutions of stochastic partial
  differential equations, Ph.D. thesis, California Institute of Technology
  (2006).

\bibitem{Hollander_etAl_2013_Nonparametric_Statistical_Methods}
M.~Hollander, D.~A. Wolfe, E.~Chicken, Nonparametric Statistical Methods, John
  Wiley \& Sons, 2013.

\bibitem{Stan_etAl_2014_CC_ConfidenceInterval}
O.~Stan, R.~Sirdey, J.~Carlier, D.~Nace, The robust binomial approach to
  chance-constrained optimization problems with application to stochastic
  partitioning of large process networks, Journal of Heuristics 20 (2014)
  1--30.

\end{thebibliography}

\end{document}